\documentclass[12pt, twoside, reqno]{amsart}
\calclayout
\usepackage[utf8]{inputenc}
\usepackage{amsmath,amssymb,amscd,amsxtra,latexsym,bbm,graphicx,epsfig,epic,
eepic,oldgerm,bm,psfrag,amsthm,mathrsfs,bm,bbm, verbatim}
\usepackage{caption}
\usepackage{subcaption}
\usepackage{tikz-cd}
\usepackage{enumitem}
\usepackage[pdftex]{hyperref}
\usepackage[noabbrev,capitalize]{cleveref}
\usetikzlibrary{decorations.pathreplacing}

\input xy
\xyoption{all}
\CompileMatrices

\theoremstyle{plain}
\newtheorem{theorem}{Theorem}[section]
\newtheorem*{theorem*}{Theorem}
\newtheorem{lemma}[theorem]{Lemma}
\newtheorem{proposition}[theorem]{Proposition}
\newtheorem{corollary}[theorem]{Corollary}
\newtheorem{definition}[theorem]{Definition}
\newtheorem*{definition*}{Definition}

\newtheorem{question}[theorem]{Question}

\newtheorem{conjecture}[theorem]{Conjecture}

\newtheoremstyle{claim}
  {3pt}
  {3pt}
  {}
  {0pt}
  {\itshape}
  {.}
  {.5em}
  {\thmname{#1}\thmnumber{ #2}\thmnote{ (#3)}}
\theoremstyle{claim}
\newtheorem*{remark*}{Remark}
\newtheorem*{remarks*}{Remarks}
\newtheorem{remark}[theorem]{Remark}

\newtheorem*{example*}{Example}
\newtheorem*{examples*}{Examples}

\theoremstyle{plain}

\newcommand{\proofend}{\hspace*{\fill} $\Box$\\}

\newcommand{\ign}[1]{}

\def\1{\:\!}
\def\2{\;\!}

\def\im{\operatorname {im}}

\def\Diffc0{\operatorname{Diff^c_0}}
\def\Symp{\operatorname{Symp}}

\def\Sympc0{\operatorname{Symp^c_0}}

\def\Int{\operatorname{int}}
\def\Ham{\operatorname{Ham}}

\def\vol{\operatorname{vol}}

\def\Aut{\operatorname{Aut}}

\def\CP{\operatorname{CP}}

\def\cd{{\mathcal D}}

\def\ch{{\mathcal H}}

\def\cl{{\mathcal L}}

\def\cR{{\mathcal R}}

\def\C{\mathbb{C}}

\def\N{\mathbb{N}}

\def\Q{\mathbb{Q}}
\def\R{\mathbb{R}}

\def\Z{\mathbb{Z}}

\def\CP{\C P}

\def\pp{\partial}

\def\ddt0{\left. \frac{d}{dt} \right\vert_{t=0}}
\def\dds0{\left. \frac{d}{ds} \right\vert_{s=0}}
\def\ddt{\frac{d}{dt} }
\def\dds{\frac{d}{ds} }

\def\Ch{\rm{Ch}}

\def\^t{^{\times}}
\def\sth{\,\vert\,}

\def\ni{\noindent}

\def\.{\mskip1mu}
\def\?{\mskip-1mu}

\def\id{\operatorname{id}}

\def\proof{\noindent {\it Proof. \;}}
\newcommand{\proofof}[1]{\ni {\it Proof of #1. }}

\begin{document}

\title[]{Lagrangian split tori in~$S^2 \times S^2$ and billiards}

\author{Jo\'e Brendel}  

\author{Joontae Kim}

\address{
D-MATH,
ETH Zürich, 
Rämistrasse 101,
8092 Zürich,
Switzerland }
\email{joe.brendel@math.ethz.ch}
\address{
Department of Mathematics and Center for Nano Materials, 35 Baekbeom-ro, Mapo-gu, Sogang University, Seoul 04107, Republic of Korea}
\email{joontae@sogang.ac.kr}

\date{\today}

\begin{abstract} 
In this paper, we classify up to Hamiltonian isotopy Lagrangian tori that split as a product of circles in $S^2 \times S^2$, when the latter is equipped with a non-monotone split symplectic form. We show that this classification is equivalent to a problem of mathematical billiards in rectangles. We give many applications, among others: (1) answering a question on Lagrangian packing numbers raised by Polterovich--Shelukhin, (2) studying the topology of the space of Lagrangian tori, and (3) determining which split tori are images under symplectic ball embeddings of Chekanov or product tori in $\R^4$.
\end{abstract}

\maketitle

\section{Introduction}

\subsection{Main result}

The problem of classifying Lagrangian submanifolds up to symplectomorphism or up to Hamiltonian diffeomorphism is very difficult in dimensions greater than two, even in simple spaces. In this paper we restrict our attention to Lagrangian tori that split as a product of circles in $S^2 \times S^2$ equipped with the symplectic form 
\begin{equation}
\label{eq:omegaalpha}
\omega_{\alpha} = 2(\alpha + 1)\omega_{S^2} \oplus 2\alpha \omega_{S^2}, 
\text{ for } \alpha > 0,
\end{equation}
where $\omega_{S^2}$ is an area form with total area $\int_{S^2} \omega_{S^2} = 1$. We denote the resulting space by $X_{\alpha} = (S^2\times S^2 , \omega_{\alpha})$. Note that up to swapping the factors and rescaling, every non-monotone symplectic form is diffeomorphic to $\omega_{\alpha}$ for some $\alpha > 0$, see \cite{LalMcD96}. \cref{rk:alphaindep} explains why we make this choice of normalization for the symplectic form. 

In $S^2(a) = (S^2,a\omega_{S^2})$, we denote by $S^1_h$, for $h \in (-\frac{a}{2},\frac{a}{2})$, the circle of constant height such that the annulus bounded by the equator and $S^1_h$ has (signed) symplectic area $h$. 

\begin{figure}
  \centering
  \begin{tikzpicture}
    \begin{scope}[scale=1,shift={(-7,0)}]
        \fill[black!5] (-2.5,-1.2) rectangle (2.5,1.2);
        \draw[black!60, thick] (-2.5,-1.2) rectangle (2.5,1.2);
        \draw[thick,black, ->] (-2.9,0)--(2.9,0);
        \draw[thick,black, ->] (0,-1.5)--(0,1.5);
        \draw[very thick,black!60, dotted] (-2,0.4)--(0,0.4);
        \draw[very thick,black!60, dotted] (-2,0)--(-2,0.4);
        \draw[very thick,black!60] (-0.1,0.4)--(0.1,0.4);
        \draw[very thick,black!60] (-2,-0.1)--(-2,0.1);
        \fill[thick, black] (-2,0.4)  circle[radius=1.8pt];
        \node at (-2,0.7) {$(x,y)$};
        \node at (2,-0.85) {$\square_{\alpha}$};
    \end{scope}

    \begin{scope}[scale=0.7,shift={(0,0)}]
        \filldraw[color=black!60, fill=black!5, very thick](0,0) circle (3);
                \begin{scope}[scale=0.765,shift={(0,-2.5)}]
                    \draw[black, very thick] (-3,0) .. controls +(0,-0.3) and +(-0.6,0).. (0,-0.4) .. controls +(0.6,0) and +(0,-0.3) .. (3,0);
			        \draw[dashed ,black, very thick] (-3,0) .. controls +(0,0.3) and +(-0.6,0).. (0,0.4) .. controls +(0.6,0) and +(0,0.3) .. (3,0);
                \end{scope}
            \draw[black!50, thick] (-3,0) .. controls +(0,-0.3) and +(-0.6,0).. (0,-0.4) .. controls +(0.6,0) and +(0,-0.3) .. (3,0);
			\draw[dashed,black!50, thick] (-3,0) .. controls +(0,0.3) and +(-0.6,0).. (0,0.4) .. controls +(0.6,0) and +(0,0.3) .. (3,0);
            \node at (2.75,-2) {$S^1_x$};
            \draw [decorate,decoration={brace,amplitude=5pt,mirror,raise=4ex}]
  (-2.2,0) -- (-2.2,-2);
            \node at (0,2) {$S^2(2\alpha + 2)$};
            \node at (-3.8,-1) {$x$};
    \end{scope}
    \node at (2.45,0) {$\times$};
    \begin{scope}[scale=0.4,shift={(10,0)}]
        \filldraw[color=black!60, fill=black!5, very thick](0,0) circle (3);
                \begin{scope}[scale=0.935,shift={(0,1)}]
                    \draw[black, very thick] (-3,0) .. controls +(0,-0.3) and +(-0.6,0).. (0,-0.4) .. controls +(0.6,0) and +(0,-0.3) .. (3,0);
			        \draw[dashed ,black, very thick] (-3,0) .. controls +(0,0.3) and +(-0.6,0).. (0,0.4) .. controls +(0.6,0) and +(0,0.3) .. (3,0);
                \end{scope}
            \draw[black!50, thick] (-3,0) .. controls +(0,-0.3) and +(-0.6,0).. (0,-0.4) .. controls +(0.6,0) and +(0,-0.3) .. (3,0);
			\draw[dashed,black!50, thick] (-3,0) .. controls +(0,0.3) and +(-0.6,0).. (0,0.4) .. controls +(0.6,0) and +(0,0.3) .. (3,0);
        \node at (-3.5,1.35) {$S^1_y$};
        \draw [decorate,decoration={brace,amplitude=5pt,mirror,raise=4ex}]
  (1.4,-0.1) -- (1.4,1);
        \node at (0,-1.75) {$S^2(2\alpha)$};
        \node at (4.1,0.45) {$y$};
    \end{scope}
  \end{tikzpicture}
  \caption{A split torus $T(x,y) \subset X_{\alpha}$ on the right-hand side and its base point in the rectangle $\square_{\alpha} = [-\alpha-1,\alpha+1] \times [-\alpha,\alpha]$ on the left-hand side. The parameters $x$ and $y$ correspond to the signed areas between the respective circles and equators.}
  \label{fig:0}
\end{figure}
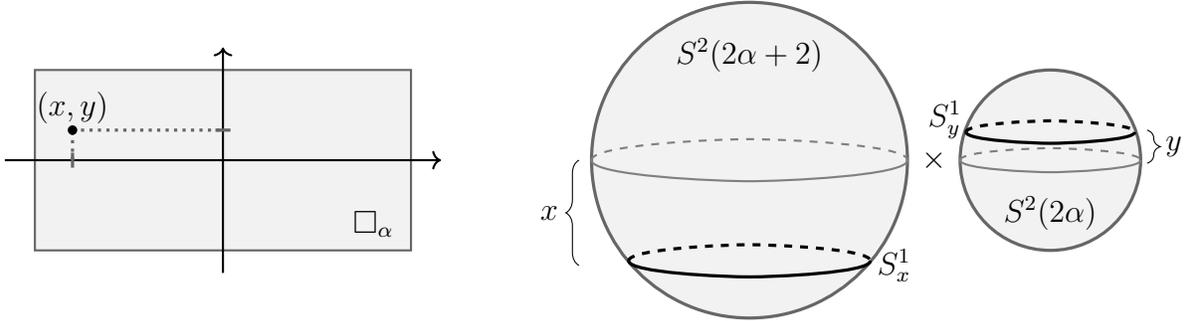

\begin{definition}
A \emph{split torus} in $X_{\alpha}$ is a torus of the form	$T(x,y) = S^1_x \times S^1_y \subset X_{\alpha}$ for $(x,y) \in (-1-\alpha, 1 +\alpha) \times (-\alpha,\alpha) \subset \R^2$.
\end{definition}

Every Lagrangian torus that splits as a product of circles is Hamiltonian isotopic to a \emph{split torus} by Hamiltonian isotopies on the $S^2$-factors. Let 
\begin{equation}
    \square_{\alpha} = [-1-\alpha,1+\alpha] \times [-\alpha,\alpha] \subset \R^2
\end{equation}
for all $\alpha > 0$. Then the set of split tori in $X_{\alpha}$ is in bijection with the interior $\Int \square_{\alpha}$, and by abuse of terminology, we will refer to \emph{a split torus whose base point is in a subset} $A \subset \Int \square_{\alpha}$ as \emph{torus in} $A$. 

The main result of this paper is a classification of split tori in $X_{\alpha}$ up to Hamiltonian diffeomorphism, which can be seen as studying the corresponding equivalence relation on points in $\Int \square_{\alpha}$. The classification up to symplectomorphism is the same, since, by connectedness of $\Symp(X_{\alpha})$ - due to Abreu \cite{Abr98} - and simple connectedness of $X_{\alpha}$, the group of symplectomorphisms agrees with the group of Hamiltonian diffeomorphisms, $\Ham (X_{\alpha})$. The analogous classification in the case of monotone $S^2 \times S^2$ was done in \cite[Proposition 5.1]{Bre23}. Throughout the paper we say that two Lagrangians $L,L'$ are \emph{equivalent} if they are Hamiltonian isotopic, and we write $L \cong L'$. If they are not Hamiltonian isotopic, we write $L \not\cong L'$.

The problem of classifying split tori in $X_{\alpha}$ is equivalent to considering certain billiard trajectories in rectangles. Let

\begin{equation}
\label{eq:Sigma}
        \Sigma = \{ x =\pm (\vert y \vert + 1) \} \cup [-1,1] \times \{0\} \subset \R^2.
\end{equation}

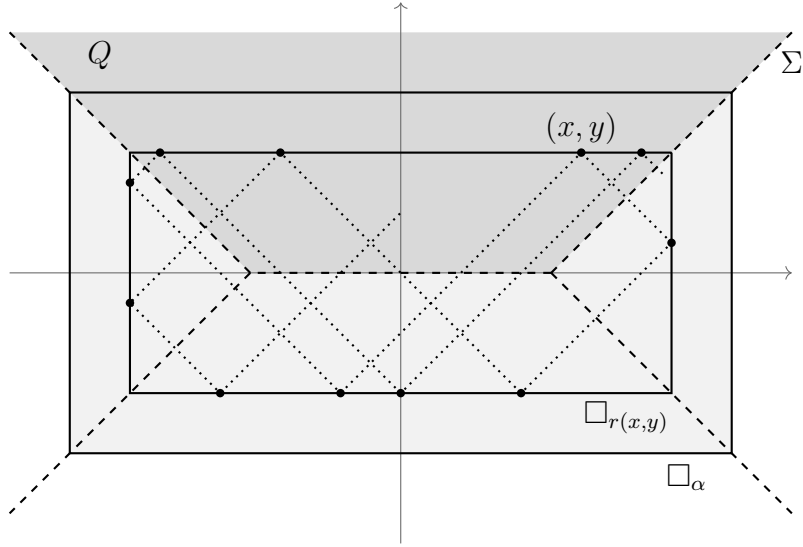
\begin{figure}
		\begin{tikzpicture}[scale=0.4]	
            \fill[black!5] (11,6)--(11,-6)--(-11,-6)--(-11,6)--(11,6);
            \fill[black!15] (5,0)--(13,8)--(-13,8)--(-5,0)--(5,0);
            \draw[thin,black!60,->] (0,-9)--(0,9);
             \draw[thin,black!60,->] (-13,0)--(13,0);
            \draw[thick,black,dashed] (5,0)--(13,8);
            \draw[thick,black,dashed] (5,0)--(13,-8);
            \draw[thick,black,dashed] (-5,0)--(-13,8);
            \draw[thick,black,dashed] (-5,0)--(-13,-8);
            \draw[thick,black,dashed] (5,0)--(-5,0);
            \draw[thick,black] (9,4)--(9,-4)--(-9,-4)--(-9,4)--(9,4);
            \draw[thick,black] (11,6)--(11,-6)--(-11,-6)--(-11,6)--(11,6);
            \draw[thick,black,dotted] (8.7,3.3)--(8,4)--(0,-4)--(-8,4)--(-9,3)--(-2,-4)--(6,4)--(9,1)--(4,-4)--(-4,4)--(-9,-1)--(-6,-4)--(0,2); 
			\fill[thick] (8,4)  circle[radius=4pt];
            \fill[thick] (0,-4)  circle[radius=4pt];
            \fill[thick] (-8,4)  circle[radius=4pt];
            \fill[thick] (-9,3)  circle[radius=4pt];
            \fill[thick] (-2,-4)  circle[radius=4pt];
            \fill[thick] (6,4)  circle[radius=4pt];
            \fill[thick] (9,1)  circle[radius=4pt];
            \fill[thick] (4,-4)  circle[radius=4pt];
            \fill[thick] (-4,4)  circle[radius=4pt];
            \fill[thick] (-9,-1)  circle[radius=4pt];
            \fill[thick] (-6,-4)  circle[radius=4pt];
            \node at (6,4.8){$(x,y)$};
            \node at (13,7){$\Sigma$};
             \node at (-10,7.25){$Q$};
            \node at (7.5,-4.8){$\square_{r(x,y)}$};
            \node at (9.5,-6.8){$\square_{\alpha}$};
		\end{tikzpicture}
	\caption{A point $(x,y) \in \square_{\alpha}$, its billiard table $\square_{r(x,y)}$, part of its good billiard trajectory (dotted line), and corresponding bouncing points (black dots). The set $\Sigma$ is indicated by dashed lines and the set $Q$ as defined in \eqref{eq:Qdefset} in grey.}
	\label{fig:1}
\end{figure}

Given $(x,y) \in \R^2$, we view the set $\square_{r(x,y)}$ as a billiard table, where $r(x,y) = \max\{\vert x \vert-1,\vert y \vert\}$. Note that $(x,y) \in \pp \square_{r(x,y)}$ and that $\pp \square_{r(x,y)} \cap \Sigma$ consists of the corners of $\square_{r(x,y)}$. See \cref{fig:1} for a sketch. The reflection law for the billiards is the usual one, where the angle of reflection is equal to the angle of incidence.

\begin{definition}
\label{def:goodbilliard}
Let $(x,y) \in \R^2$. If $(x,y) \notin \Sigma$, the \emph{good billiard trajectory of $(x,y)$} is the unique billiard trajectory in $\square_{r(x,y)}$ that
\begin{itemize}
\item[{\rm (1)}] has incidence angle $\frac{\pi}{4}$ with $\pp \square_{r(x,y)}$,
\item[{\rm (2)}] contains $(x,y)$ as bouncing point, 
\item[{\rm (3)}] is reflected backwards in the same direction when it hits a corner of the billiard table $\pp \square_{r(x,y)}$. 
\end{itemize}
In that case, we call a bouncing point of the good billiard trajectory of $(x,y)$ \emph{admissible} if it is not in $\Sigma$, i.e.\ if it is not a corner point of the corresponding billiard table. If $(x,y) \in \Sigma$ we define the \emph{good billiard trajectory of $(x,y)$} as stationary and call its unique bouncing point \emph{admissible}.
\end{definition}

This terminology is chosen so that our results can be formulated in a unified way. See \cref{fig:1} for a sketch of a good billiard trajectory. Note that good billiard trajectories can have very few admissible bouncing points, even if $(x,y) \notin \Sigma$. For example, the good billiard trajectory of $(x,y)=(0,1)$ (in case $\alpha > 1$) has a unique admissible bouncing point.

\begin{theorem}
\label{thm:main}
Let $T(x,y),T(x',y') \subset X_{\alpha}$ be split tori. Then $T(x,y) \cong T(x',y')$ if and only if the good billiard trajectory emanating from $(x,y)$ has any of the points $(\pm x', \pm y')$ as admissible bouncing point.
\end{theorem}

This result is proved in \cref{sec:mainthm} along with its refinement \cref{thm:main2}. See \S\ref{ssec:introspaceoflags}-\S\ref{ssec:introballemb} for various applications of this theorem. 

\begin{remark}
\label{rk:alphaindep}
Strikingly, this implies that the classification of split tori is independent of $\alpha > 0$. Indeed, let $(x,y) \in \R^2$ and let $\alpha, \alpha' > r(x,y)$. Since the statement of \cref{thm:main} depends on the billiard table $\square_{r(x,y)}$ only, two split tori $T(x,y),T(x',y') \subset X_{\alpha}$ are equivalent if and only if their counterparts $T(x,y),T(x',y') \subset X_{\alpha'}$ are equivalent. This is a by-product of the choice of symplectic form in \eqref{eq:omegaalpha} and explains this convention. 
\end{remark}

\begin{remark}
For every $\alpha > 0$, there are Lagrangian tori in $X_{\alpha}$ that cannot be mapped to a split torus by a symplectomorphism. This follows for example from~\cite[Theorem C]{BreHauSch23}. Therefore, this classification result concerns a strict subset of all Lagrangian tori in $X_{\alpha}$.
\end{remark}

\begin{remark}
\label{rk:Hirzebruch}
Split tori can be viewed as toric fibres of the standard moment map on $X_{\alpha}$, which has $\square_{\alpha}$ as image. For some $\alpha > 0$, certain even Hirzebruch surfaces yield other toric structures on $X_{\alpha}$, see for example \cite{Kar03}, and it is natural to ask for a classification of their toric fibres up to Hamiltonian diffeomorphism. Through the use of almost toric fibrations, these toric structures are all related by mutation, see \cite[\S 6.2]{Sym03}. Although this requires some additional work which is beyond the scope of the present paper, this approach proves that the set of toric fibres of even Hirzebruch surfaces coincides with the set of split tori up to Hamiltonian diffeomorphism\footnote{In the sense that for every toric fibre of an even Hirzebruch surface, there is a Hamiltonian diffeomorphism which maps it to a split torus. The mutation of the base diagram allows one to identify which split torus it is.}. Therefore, our classification yields a classification of toric fibres in these alternative toric structures, too. Similarly, to achieve a classification of all toric fibres of odd Hirzebruch surfaces, it would be sufficient to classify toric fibres in the one-fold blow-up of $\CP^2$ for all possible symplectic forms; see \cref{q:class1}.
\end{remark}

\begin{remark}
Recall the elementary fact that if the ratio of width and height of a rectangular billiard table is irrational, then the set of bouncing points of any good billiard trajectory is dense in the boundary of that billiard table. By \cref{thm:main}, this implies that there are sets of equivalent tori in $X_{\alpha}$ which are dense in the boundary of certain rectangles in $\square_{\alpha}$. This will be used repeatedly in applications, for example to produce infinite packings by Lagrangian tori of $X_{\alpha}$ and to exhibit various pathologies of the space of Lagrangian tori.
\end{remark}

\subsection{Hamiltonian monodromy}
A question closely related to the classification of split tori is determining their Hamiltonian monodromy group. 

\begin{definition}
\label{def:hammonodromygroup}
Let~$L \subset X$ be a Lagrangian in a symplectic manifold $(X,\omega)$. The group
\begin{equation}
    \label{eq:hammonodromygroup}
	\ch_L = \left\{ (\psi\vert_L)_* \in \Aut(H_1(L)) \, \vert \, \psi \in \Ham(X,\omega), \psi(L) = L \right\}
\end{equation} 
is called the \emph{Hamiltonian monodromy group of~$L$}.
\end{definition}

The Hamiltonian monodromy group has been studied quite extensively, see \cite{Che96, Yau09, HuLalLec11, Ono15, AugSmiWor22}. For a thorough discussion of the existing results, we refer to the introduction of \cite{AugSmiWor22}. In \cite{Bre23}, the first named author has determined the Hamiltonian monodromy groups of all toric fibres in certain simple four-dimensional examples, such as monotone $S^2 \times S^2$. This extends the earlier result by Ono \cite{Ono15} for the Clifford torus in monotone $S^2 \times S^2$ to non-monotone split tori. 

In this paper we study the Hamiltonian monodromy group of split tori in the spaces $X_{\alpha}$ and determine it completely in many cases, though not in all. We find that, in a suitable basis of its first homology, every element in the Hamiltonian monodromy group of a split torus in $X_{\alpha}$ is of the form
\begin{equation}
    \begin{pmatrix}
		\pm 1 & 0 \\
		2k & \pm 1
	\end{pmatrix}
\end{equation}
for some $k \in \Z$. The Hamiltonian monodromy group subtly depends on arithmetic properties of the base point $(x,y)$ of $L = T(x,y)$, and many different isomorphism types of groups arise in this problem; for example $\{1\}$, $\Z_2$, but also infinite groups such as $\Z$ and various semi-direct products of $\Z_2$ with $\Z$. Let us point out that in the case of split tori in monotone $S^2 \times S^2$, only finite groups appear. Our full result is somewhat cumbersome to state; see \cref{thm:monodromy} for more details.

\subsection{Methods of proof}

Split tori are toric fibres of the standard toric structure $X_{\alpha} \rightarrow \square_{\alpha}$ and we make extensive use of toric geometry. In particular, based on~\cite{Bre23}, we use symmetric probes to construct Hamiltonian diffeomorphisms realizing equivalences of fibres and given elements of the Hamiltonian monodromy group. Probes were introduced in~\cite{McD11} and symmetric probes in~\cite{AbrBorMcD14} for the purpose of displacing toric fibres. The obstructions are derived from Chekanov's obstructions for product tori in $\R^{2n}$ to be symplectomorphic via a lifting trick from toric geometry; see \cite[Section 4]{Bre23}. A key novelty in this paper is the use of properties of billiards to take care of the somewhat involved combinatorics arising in the context of non-monotone $S^2 \times S^2$. A key insight is that, up to reflections $(x,y) \mapsto (\pm x, \pm y)$, only iterated symmetric probes realizing billiard trajectories are needed to complete the classification. In \cite{Bre23}, the first named author has conjectured that two toric fibres in a closed toric symplectic manifold are Hamiltonian isotopic if and only if there is a sequence of symmetric probes between them. \cref{thm:main} confirms this conjecture in the case of non-monotone $S^2 \times S^2$. 

\begin{corollary}
    Conjecture 1.3 in \cite{Bre23} holds for $S^2 \times S^2$ equipped with a non-monotone split symplectic form.
\end{corollary}

Similarly, we suspect that every non-trivial element in the Hamiltonian monodromy group of a split torus is generated by a sequence of symmetric probes starting and ending on the same split torus. Note however that for \eqref{eq:mon4.5} and \eqref{eq:mon8}, \cref{thm:monodromy} does not completely determine the monodromy group. Therefore this problem is still open for these cases.

\subsection{Application: Space of Lagrangians}
\label{ssec:introspaceoflags}
Let $L=T(x,y)$ with $r(x,y) \in \R \setminus \Q$ and $(x,y) \notin \Sigma$. Since the ratio of width to height of the corresponding billiard table $\square_{r(x,y)}$ is irrational, there is a dense subset of points in $\pp \square_{r(x,y)}$ such that the corresponding tori are Hamiltonian isotopic. Inspired by the recent work of Chassé--Leclercq \cite{ChaLec24}, we deduce the following two immediate consequences of \cref{thm:main}, which give counterexamples to Conjectures A and B in \cite{ChaLec24} in dimension four.

\begin{corollary}
    \label{cor:weinstein}
    In every $X_{\alpha}$, there is a split torus $L \subset X_{\alpha}$ such that every Weinstein neighbourhood of $L$ contains a (and thus infinitely many) split torus $L'$ with $L \cong L'$ and $L \cap L' = \varnothing$.
\end{corollary}

By the density of bouncing points in $\pp \square_{r(x,y)}$ with $r(x,y) \in \R \setminus \Q$, we can select a sequence of points which converges to a point not in the set of bouncing points to obtain the following.

\begin{corollary}
    \label{cor:convergence}
    In every $X_{\alpha}$, there is a sequence of split tori $L_i \subset X_{\alpha}$ which are mutually Hamiltonian isotopic, $L_i \cong L_j$, and which converge in the $C^{\infty}$-topology to a split torus $L \subset X_{\alpha}$ with $L \not\cong L_i$. 
\end{corollary}

In particular, the orbit under $\Ham(X_{\alpha})$ of such a torus in the space of Lagrangians is neither $C^{\infty}$-closed, nor locally path connected. Both these phenomena occur in any symplectic manifold of dimension greater than or equal to six. As was proved in \cite{ChaLec24}, this follows from Chekanov's classification \cite{Che96} of product tori in $\R^{2n}$, see also \cite{Bre23} for a slightly different point of view. 

\begin{remark}
Note that the Lebesgue measure on $\square_{\alpha} \subset \R^2$ induces a natural measure on the set of split tori in $X_{\alpha}$. The subset of tori having the properties from the Corollaries \ref{cor:weinstein} and \ref{cor:convergence} have full measure. 
\end{remark}

\begin{remark}
As was pointed out to us by Felix Schlenk, one can also use \cite[Theorem~1.5]{CheSch16} to produce results of this type in dimension four. Note however, that their methods work only for small tori contained in a Darboux chart, and, in the special case of $X_{\alpha}$, only for certain $\alpha > 0$.
\end{remark}

\subsection{Application: Lagrangian packing}

\label{ssec:introLagpacking}

\begin{definition}
Let $L$ be a Lagrangian in a symplectic manifold $(X,\omega)$. The \emph{packing number} of $L$ is defined as 
\begin{equation}
	\#_P(L) = \sup\{ N \in \N \sth \exists \, \phi_1,\ldots,\phi_N \in \Ham(X,\omega),\, \phi_i(L) \cap \phi_j(L) = \varnothing,\, \forall i \neq j \}.
\end{equation}
\end{definition}

The following remarkable result is due to Polterovich--Shelukhin, see \cite[Theorem C]{PolShe23}. It uses ideas introduced by Mak--Smith \cite{MakSmi21}, who proved the non-displaceability of certain links of Lagrangian tori in $X_{\alpha}$ by computing Floer homology in the symmetric product.

\begin{theorem}[Polterovich--Shelukhin]
\label{thm:MSPS}
Let $\alpha > 0$, $L = T(x,0) \in X_{\alpha}$ and $k \in \N_{\geqslant 2}$. Then 
\begin{equation}
    \label{eq:PSmain}
    \vert x \vert < \frac{k-1}{k+1} \quad \Rightarrow \quad \#_P(L) \leqslant k.
\end{equation}
\end{theorem}

In \S\ref{ssec:discussionMSPS}, we give details on how to adapt \cite[Theorem C]{PolShe23} to our notational conventions. For any $\vert x \vert \in \left( \frac{k-2}{k} ,  \frac{k-1}{k+1} \right)$, the bound \eqref{eq:PSmain} is optimal for small enough $\alpha > 0$. Indeed, in that case, a $k$-packing can be realized by packing the first factor in $X_{\alpha} = S^2 \times S^2$ by $k$ circles. 

Since $\frac{k-1}{k+1} \rightarrow 1$ for $k \rightarrow + \infty$, \cref{thm:MSPS} implies that the packing number of every split torus $T(x,0)$ with $\vert x \vert < 1$ is finite. We show that this rigidity breaks down for $\vert x \vert > 1$, answering \cite[Question 18]{PolShe23} and showing that, up to the critical case $\vert x \vert = 1$, which remains open, the results in \cite{PolShe23} are sharp. Indeed, since distinct split tori are disjoint, \cref{thm:main} yields many interesting packings by split tori. In particular, we find $\#_P(L) = \infty$ for split tori $L = T(x,0)$ with $x \in (1,1+\alpha) \setminus \Q$. This is a special case of the following, which itself follows immediately from \cref{thm:main} and the fact that good billiard trajectories have infinitely many bouncing points whenever the ratio of width and height of the billiard table is irrational.

\begin{corollary}
\label{cor:infpacking}
Let $L = T(x,y) \subset X_{\alpha}$ be a split torus with $r(x,y) \in \R \setminus \Q$. Then 
\begin{equation}
    \#_P(L) = \infty.
\end{equation}
\end{corollary}

Note, again, that this shows that generic split tori have this property. As a consequence of \cref{thm:main}, we determine all packing numbers of $X_{\alpha}$ \emph{by disjoint split tori}, see \cref{cor:lowerbound} and \cref{prop:toricpacking}. This gives a lower bound on the actual Lagrangian packing number, which is optimal in some cases but not in others, see also the discussion in \S\ref{ssec:qpacking}. Note that even in the case $r(x,y) \in \Q$, this lower bound tends to be much larger than what one would obtain by a naive packing of circles into one of the sphere factors. There are many open questions surrounding this set of ideas, see Questions \ref{q:packing1},\ref{q:packing2}.

\subsection{Application: Late Poincaré recurrence}

\label{ssec:introlatereturn}

The Lagrangian Poincaré recurrence conjecture can be stated as follows.

\begin{conjecture}[Lagrangian Poincaré recurrence Conjecture]
Let $\psi \in \Ham(X,\omega)$ and let $L \subset (X,\omega)$ be a closed Lagrangian submanifold. Then there exists a sequence $k_i \rightarrow \infty$ of natural numbers such that $\psi^{k_i}(L) \cap L \neq \varnothing$. Furthermore, the density of the sequence $k_i$ is related to a symplectic capacity of $L$.
\end{conjecture}

We refer to \cite[Conjecture 5.6]{GinGur18}, where the conjecture is attributed to Ginzburg and Viterbo independently. Recently, Broćić--Shelukhin \cite{BroShe24} have found counterexamples to the LPR Conjecture in all symplectic manifolds of dimensions greater than or equal to six and Schmitz \cite{Sch24} in some symplectic manifolds of dimension four. None of these results apply to the case of $S^2 \times S^2$. Using our methods, we prove a related but weaker result: For certain fixed split tori $L \subset X_{\alpha}$, we find Hamiltonian diffeomorphisms realizing arbitrarily large (but finite) return time.

\begin{proposition}
\label{prop:latereturnweak}
    Let $\alpha > 0$ and let $L = T(x,y) \subset X_{\alpha}$ be a split torus with $r(x,y) \in \R \setminus \Q$. For every $N \in \N$, there is $\psi \in \Ham(X_{\alpha})$ with 
    \begin{equation}
        \psi^k(L) \cap L = \varnothing, \quad
        \text{for all } k \in \{1,\ldots,N\}.
    \end{equation}
\end{proposition}

This follows straightforwardly from iterated symmetric probes, see the proof in \cref{sec:latereturn}. With a bit more work, we can prove a stronger version. For $\phi \in \Ham(X_{\alpha})$ and $N \in \N$, let 
\begin{equation}
    \label{eq:mathcaltdef}
	\mathcal{T}_{\phi,N} = \{ (x,y) \in \Int \square_{\alpha} \, \vert \, \phi^k(T(x,y)) \cap T(x,y) = \varnothing,\, \forall k \in \{1,\ldots, N\}\},
\end{equation}
i.e.\ this is the set of all split tori whose images under $\phi$ are disjoint under the $N$-th iterates. We define a family $\{\psi_{\delta}\}_{\delta > 0} \subset \Ham (X_{\alpha})$ out of which one can pick an element making this set arbitrarily large for arbitrarily large $N$. By $\lambda$ we denote Lebesgue measure on $\square_{\alpha} \subset \R^2$ with normalization $\lambda(\square_{\alpha})=1$.

\begin{theorem}
\label{thm:latereturn}
For all $\varepsilon > 0$ and $N \in \N$, there is $\delta >0$ such that 
\begin{equation}
	\lambda(\mathcal{T}_{\psi_{\delta},N}) > 1 - \varepsilon.
\end{equation}
\end{theorem}

\subsection{Application: Lagrangian tori from symplectic ball embeddings}

\label{ssec:introballemb}

For any Lagrangian torus $L \subset X$, one can ask whether $L$ is the image of a Lagrangian torus under a symplectic ball embedding and, if so, which Hamiltonian isotopy type of torus in the ball it corresponds to. 

\begin{definition}
We call a Lagrangian torus $L \subset (X,\omega)$ \emph{ball-embeddable} if there is a symplectic ball embedding $\varphi \colon B^{2n}(r) \hookrightarrow X$ with $L \subset \im \varphi$.
\end{definition}

By $B^{2n}(r) = \{ x_1^2 + y_1^2 + \ldots + x_n^2 + y_n^2 \leqslant \pi r^2 \} \subset (\R^{2n}, \omega_0)$, we denote the \emph{closed} $2n$-ball of capacity $r$. In dimension four, we consider three types of Lagrangian tori in $\R^4$, namely \emph{nonmonotone product tori}, i.e.\ tori of the form $\Theta(a_1,a_2) = S^1(a_1) \times S^1(a_2)$ with $a_1 \neq a_2$, \emph{Clifford tori}, i.e.\ product tori of the form $T_{\rm Cl}(a) = \Theta(a,a)$, and \emph{Chekanov tori}, $T_{\rm Ch}(a)$, for which we give a definition in \cref{sec:ballembeddable}. 

\begin{definition}
\label{def:types}
Let $L \subset X^4$ be a ball-embeddable torus. We say that $L$ is \emph{ball-Clifford} if it is contained in the image of a symplectic ball embedding $\varphi$ into $X$ such that $\varphi^{-1}(L) \subset \R^4$ is a Clifford torus. The terms \emph{ball-Chekanov} and \emph{ball-nonmonotone} are defined analogously.
\end{definition}

A priori, a ball-embeddable torus may be of none of these three types; or of several types simultaneously. In fact, we will see that the latter occurs in the case of $X = X_{\alpha}$. We first note the following.

\begin{theorem}
\label{thm:notexotic}
For every symplectic ball embedding $\varphi \colon B^4(r) \hookrightarrow X_{\alpha}$ and $a < \frac{r}{2}$, we have 
\begin{equation}
    \label{eq:notexotic}
	\varphi(T_{\Ch}(a)) \cong T(\alpha - a -1, \alpha - a).
\end{equation}
\end{theorem}

In other words, \emph{ball-Chekanov tori in $X_{\alpha}$ are not exotic}. Note that this stands in contrast with the case of monotone $S^2 \times S^2$, where ball-Chekanov tori \emph{are} exotic, i.e.\ not Hamiltonian isotopic to a split torus. For the monotone Chekanov torus (though it is not ball embeddable), this was proved by Chekanov--Schlenk \cite{CheSch10}. For the those non-monotone Chekanov tori that are ball-embeddable, this is a recent result by Lou \cite{Lou24}.

Our main result about ball-embeddability is completely determining which split tori are ball-Clifford, which are ball-Chekanov and which are ball-nonmonotone.

\begin{theorem}
\label{thm:ballcases}
A split torus $L = T(x,y) \subset X_{\alpha}$ is 
\begin{enumerate}
    \item ball-Clifford if and only if\footnote{Recall the definition of $\Sigma$ from \eqref{eq:Sigma}, see also \cref{fig:1}.} $(x,y) \in \Sigma$ with $y \neq 0$;
    \item ball-Chekanov if and only if $(x,y) \notin \Sigma$ and the good billiard trajectory of $(x,y)$ contains a corner of the billiard table $\square_{r(x,y)}$ as bouncing point; 
    \item ball-nonmonotone if and only if $(x,y)$ is contained in the complement of the set $\Sigma \cup \cd$, where
        \begin{equation}
            \label{eq:nonmonotonecomplement}
            \cd = \left\{ \left. \left( \frac{k-2l-1}{k}, \pm \frac{1}{k} \right) \right\vert k \in \N_{\geqslant 1}, l \in \{0,\ldots,k-1\} \right\}.
        \end{equation}
\end{enumerate}
\end{theorem}

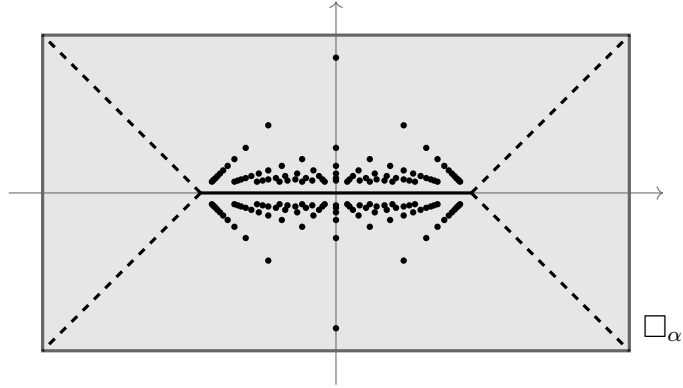
\begin{figure}
		\begin{tikzpicture}[scale=0.3]	
            \fill[black!10] (-13,-7) rectangle (13,7);
            \draw[black!60, very thick] (-13,-7) rectangle (13,7);
            \draw[thin,black!60,->] (0,-8.5)--(0,8.5);
             \draw[thin,black!60,->] (-14.5,0)--(14.5,0);
            \draw[very thick,black,dashed] (6,0)--(13,7);
            \draw[very thick,black,dashed] (6,0)--(13,-7);
            \draw[very thick,black,dashed] (-6,0)--(-13,7);
            \draw[very thick,black,dashed] (-6,0)--(-13,-7);
            \draw[very thick,black] (6,0)--(-6,0);
			
            \fill[thick] (0,6)  circle[radius=4pt];
            
            \fill[thick] (-3,3)  circle[radius=4pt];
            \fill[thick] (3,3)  circle[radius=4pt];
            
            \fill[thick] (-4,2)  circle[radius=4pt];
            \fill[thick] (0,2)  circle[radius=4pt];
            \fill[thick] (4,2)  circle[radius=4pt];
            
            \fill[thick] (-4.5,1.5)  circle[radius=4pt];
            \fill[thick] (-1.5,1.5)  circle[radius=4pt];
            \fill[thick] (1.5,1.5)  circle[radius=4pt];
            \fill[thick] (4.5,1.5)  circle[radius=4pt];
            
            \fill[thick] (-4.8,1.2)  circle[radius=4pt];
            \fill[thick] (-2.4,1.2)  circle[radius=4pt];
            \fill[thick] (0,1.2)  circle[radius=4pt];
            \fill[thick] (2.4,1.2)  circle[radius=4pt];
            \fill[thick] (4.8,1.2)  circle[radius=4pt];
            
            \fill[thick] (-5,1)  circle[radius=4pt];
            \fill[thick] (-3,1)  circle[radius=4pt];
            \fill[thick] (-1,1)  circle[radius=4pt];
            \fill[thick] (1,1)  circle[radius=4pt];
            \fill[thick] (3,1)  circle[radius=4pt];
            \fill[thick] (5,1)  circle[radius=4pt];
            
            \fill[thick] (-5.14,0.86)  circle[radius=4pt];
            \fill[thick] (-3.42,0.86)  circle[radius=4pt];
            \fill[thick] (-1.72,0.86)  circle[radius=4pt];
            \fill[thick] (-0,0.86)  circle[radius=4pt];
            \fill[thick] (1.72,0.86)  circle[radius=4pt];
            \fill[thick] (3.42,0.86)  circle[radius=4pt];
            \fill[thick] (5.14,0.86)  circle[radius=4pt];
            
            \fill[thick] (-5.25,0.75)  circle[radius=4pt];
            \fill[thick] (-3.75,0.75)  circle[radius=4pt];
            \fill[thick] (-2.25,0.75)  circle[radius=4pt];
            \fill[thick] (-0.75,0.75)  circle[radius=4pt];
             \fill[thick] (5.25,0.75)  circle[radius=4pt];
            \fill[thick] (3.75,0.75)  circle[radius=4pt];
            \fill[thick] (2.25,0.75)  circle[radius=4pt];
            \fill[thick] (0.75,0.75)  circle[radius=4pt];

            \fill[thick] (-5.33,0.66)  circle[radius=4pt];
            \fill[thick] (-4,0.66)  circle[radius=4pt];
            \fill[thick] (-2.66,0.66)  circle[radius=4pt];
            \fill[thick] (-1.33,0.66)  circle[radius=4pt];
            \fill[thick] (0,0.66)  circle[radius=4pt];
            \fill[thick] (5.33,0.66)  circle[radius=4pt];
            \fill[thick] (4,0.66)  circle[radius=4pt];
            \fill[thick] (2.66,0.66)  circle[radius=4pt];
            \fill[thick] (1.33,0.66)  circle[radius=4pt];

            \fill[thick] (-5.4,0.6)  circle[radius=4pt];
            \fill[thick] (-4.2,0.6)  circle[radius=4pt];
            \fill[thick] (-3,0.6)  circle[radius=4pt];
            \fill[thick] (-1.8,0.6)  circle[radius=4pt];
            \fill[thick] (-0.6,0.6)  circle[radius=4pt];
            \fill[thick] (5.4,0.6)  circle[radius=4pt];
            \fill[thick] (4.2,0.6)  circle[radius=4pt];
            \fill[thick] (3,0.6)  circle[radius=4pt];
            \fill[thick] (1.8,0.6)  circle[radius=4pt];
            \fill[thick] (0.6,0.6)  circle[radius=4pt];

            \fill[thick] (-5.45,0.55)  circle[radius=4pt];
            \fill[thick] (-4.35,0.55)  circle[radius=4pt];
            \fill[thick] (-3.25,0.55)  circle[radius=4pt];
            \fill[thick] (-2.15,0.55)  circle[radius=4pt];
            \fill[thick] (-1.05,0.55)  circle[radius=4pt];
            \fill[thick] (0,0.55)  circle[radius=4pt];
            \fill[thick] (5.45,0.55)  circle[radius=4pt];
            \fill[thick] (4.35,0.55)  circle[radius=4pt];
            \fill[thick] (3.25,0.55)  circle[radius=4pt];
            \fill[thick] (2.15,0.55)  circle[radius=4pt];
            \fill[thick] (1.05,0.55)  circle[radius=4pt];

            \fill[thick] (-5.5,0.5)  circle[radius=4pt];
            \fill[thick] (-4.5,0.5)  circle[radius=4pt];
            \fill[thick] (-3.5,0.5)  circle[radius=4pt];
            \fill[thick] (-2.5,0.5)  circle[radius=4pt];
            \fill[thick] (-1.5,0.5)  circle[radius=4pt];
            \fill[thick] (-0.5,0.5)  circle[radius=4pt];
            \fill[thick] (5.5,0.5)  circle[radius=4pt];
            \fill[thick] (4.5,0.5)  circle[radius=4pt];
            \fill[thick] (3.5,0.5)  circle[radius=4pt];
            \fill[thick] (2.5,0.5)  circle[radius=4pt];
            \fill[thick] (1.5,0.5)  circle[radius=4pt];
            \fill[thick] (0.5,0.5)  circle[radius=4pt];

            \fill[thick] (0,-6)  circle[radius=4pt];
            
            \fill[thick] (-3,-3)  circle[radius=4pt];
            \fill[thick] (3,-3)  circle[radius=4pt];
            
            \fill[thick] (-4,-2)  circle[radius=4pt];
            \fill[thick] (0,-2)  circle[radius=4pt];
            \fill[thick] (4,-2)  circle[radius=4pt];
            
            \fill[thick] (-4.5,-1.5)  circle[radius=4pt];
            \fill[thick] (-1.5,-1.5)  circle[radius=4pt];
            \fill[thick] (1.5,-1.5)  circle[radius=4pt];
            \fill[thick] (4.5,-1.5)  circle[radius=4pt];
            
            \fill[thick] (-4.8,-1.2)  circle[radius=4pt];
            \fill[thick] (-2.4,-1.2)  circle[radius=4pt];
            \fill[thick] (0,-1.2)  circle[radius=4pt];
            \fill[thick] (2.4,-1.2)  circle[radius=4pt];
            \fill[thick] (4.8,-1.2)  circle[radius=4pt];
            
            \fill[thick] (-5,-1)  circle[radius=4pt];
            \fill[thick] (-3,-1)  circle[radius=4pt];
            \fill[thick] (-1,-1)  circle[radius=4pt];
            \fill[thick] (1,-1)  circle[radius=4pt];
            \fill[thick] (3,-1)  circle[radius=4pt];
            \fill[thick] (5,-1)  circle[radius=4pt];
            
            \fill[thick] (-5.14,-0.86)  circle[radius=4pt];
            \fill[thick] (-3.42,-0.86)  circle[radius=4pt];
            \fill[thick] (-1.72,-0.86)  circle[radius=4pt];
            \fill[thick] (-0,-0.86)  circle[radius=4pt];
            \fill[thick] (1.72,-0.86)  circle[radius=4pt];
            \fill[thick] (3.42,-0.86)  circle[radius=4pt];
            \fill[thick] (5.14,-0.86)  circle[radius=4pt];
            
            \fill[thick] (-5.25,-0.75)  circle[radius=4pt];
            \fill[thick] (-3.75,-0.75)  circle[radius=4pt];
            \fill[thick] (-2.25,-0.75)  circle[radius=4pt];
            \fill[thick] (-0.75,-0.75)  circle[radius=4pt];
             \fill[thick] (5.25,-0.75)  circle[radius=4pt];
            \fill[thick] (3.75,-0.75)  circle[radius=4pt];
            \fill[thick] (2.25,-0.75)  circle[radius=4pt];
            \fill[thick] (0.75,-0.75)  circle[radius=4pt];

            \fill[thick] (-5.33,-0.66)  circle[radius=4pt];
            \fill[thick] (-4,-0.66)  circle[radius=4pt];
            \fill[thick] (-2.66,-0.66)  circle[radius=4pt];
            \fill[thick] (-1.33,-0.66)  circle[radius=4pt];
            \fill[thick] (0,-0.66)  circle[radius=4pt];
            \fill[thick] (5.33,-0.66)  circle[radius=4pt];
            \fill[thick] (4,-0.66)  circle[radius=4pt];
            \fill[thick] (2.66,-0.66)  circle[radius=4pt];
            \fill[thick] (1.33,-0.66)  circle[radius=4pt];

            \fill[thick] (-5.4,-0.6)  circle[radius=4pt];
            \fill[thick] (-4.2,-0.6)  circle[radius=4pt];
            \fill[thick] (-3,-0.6)  circle[radius=4pt];
            \fill[thick] (-1.8,-0.6)  circle[radius=4pt];
            \fill[thick] (-0.6,-0.6)  circle[radius=4pt];
            \fill[thick] (5.4,-0.6)  circle[radius=4pt];
            \fill[thick] (4.2,-0.6)  circle[radius=4pt];
            \fill[thick] (3,-0.6)  circle[radius=4pt];
            \fill[thick] (1.8,-0.6)  circle[radius=4pt];
            \fill[thick] (0.6,-0.6)  circle[radius=4pt];

            \fill[thick] (-5.45,-0.55)  circle[radius=4pt];
            \fill[thick] (-4.35,-0.55)  circle[radius=4pt];
            \fill[thick] (-3.25,-0.55)  circle[radius=4pt];
            \fill[thick] (-2.15,-0.55)  circle[radius=4pt];
            \fill[thick] (-1.05,-0.55)  circle[radius=4pt];
            \fill[thick] (0,-0.55)  circle[radius=4pt];
            \fill[thick] (5.45,-0.55)  circle[radius=4pt];
            \fill[thick] (4.35,-0.55)  circle[radius=4pt];
            \fill[thick] (3.25,-0.55)  circle[radius=4pt];
            \fill[thick] (2.15,-0.55)  circle[radius=4pt];
            \fill[thick] (1.05,-0.55)  circle[radius=4pt];

            \fill[thick] (-5.5,-0.5)  circle[radius=4pt];
            \fill[thick] (-4.5,-0.5)  circle[radius=4pt];
            \fill[thick] (-3.5,-0.5)  circle[radius=4pt];
            \fill[thick] (-2.5,-0.5)  circle[radius=4pt];
            \fill[thick] (-1.5,-0.5)  circle[radius=4pt];
            \fill[thick] (-0.5,-0.5)  circle[radius=4pt];
            \fill[thick] (5.5,-0.5)  circle[radius=4pt];
            \fill[thick] (4.5,-0.5)  circle[radius=4pt];
            \fill[thick] (3.5,-0.5)  circle[radius=4pt];
            \fill[thick] (2.5,-0.5)  circle[radius=4pt];
            \fill[thick] (1.5,-0.5)  circle[radius=4pt];
            \fill[thick] (0.5,-0.5)  circle[radius=4pt];

            \node at (14.5,-6){$\square_{\alpha}$};
		\end{tikzpicture}
	\caption{Illustration of \cref{thm:ballcases}. The set of ball-Clifford split tori is indicated by dashed lines. The dots represent the second set in the union \eqref{eq:nonmonotonecomplement} i.e.\ the tori which are ball-Chekanov but not ball-nonmonotone. The central segment drawn by a solid line consists of the points which are of neither of the three types in \cref{def:types}. The set of ball-nonmonotone tori is the complement of the union of dashed lines, the central segment and the dots. The set of ball-Chekanov tori is dense in $\square_{\alpha}$ and hence cannot be meaningfully represented.}
	\label{fig:1b}
\end{figure}

See \cref{fig:1b} for an illustration. The proof of \cref{thm:ballcases} heavily relies on the fact that the space of symplectic ball embeddings into $X_{\alpha}$ is connected, which goes back to McDuff~\cite{McD98}; see also \cref{lem:ballconn}.

\cref{thm:ballcases} has interesting consequences. First note that, among split tori, being ball-Chekanov and ball-Clifford are mutually exclusive. The same holds for ball-nonmonotone and ball-Clifford. However the sets of \emph{ball-Chekanov} tori and \emph{ball-nonmonotone} tori intersect non-trivially. More precisely, we find the following. 

\begin{corollary}
\label{cor:cheknonmonotone}
Let $\varphi \colon B^4(r) \hookrightarrow X_{\alpha}$ be a symplectic ball embedding and $T_{\rm Ch}(a) \subset B^4(r)$. There is a product torus $\Theta(b,c) \subset B^4(r)$ with
\begin{equation}
    \label{eq:chekequiv}
	\varphi(T_{\Ch}(a)) \cong \varphi(\Theta(b,c)),
\end{equation}
if and only if $a \notin \left\{ \alpha - \frac{1}{k} \sth k \in \N_{\geqslant 1} \right\}$. In this case, the torus $\Theta(b,c)$ is nonmonotone, i.e. $b \neq c$.
\end{corollary}

Since the Chekanov torus is exotic in $\R^4$, an isotopy as in \eqref{eq:chekequiv} cannot be supported in any extension of the symplectic ball embedding to a larger domain in $\R^4$. Hence any such Hamiltonian isotopy in $X_{\alpha}$ has to intersect a divisor at infinity. Another consequence of \cref{thm:ballcases} is that the split tori which are of neither of the three types in \cref{def:types} are exactly those over the central segment $[-1,1] \times \{0\}$; see \cref{fig:1b}. We do not know, in general, if these tori are ball-embeddable or not. For a discussion of this and other open questions surrounding Lagrangian tori coming from ball embeddings, we refer to \S\ref{ssec:qballembeddable}.

\subsection{Structure of the paper} The heart of this paper is \cref{sec:mainthm}, where we review symmetric probes, introduce methods from billiards, and prove our main result, \cref{thm:main}. In \cref{sec:monodromy} we explore the relationship between billiards and the Hamiltonian monodromy group of split tori, prove \cref{thm:monodromy} and discuss why our methods are insufficient in some special cases. In \cref{sec:ballembeddable}, we prove the results about ball-embeddable tori outlined in \S\ref{ssec:introballemb}. In \cref{sec:Lagpacking}, we prove the results in \S\ref{ssec:introLagpacking} and in \cref{sec:latereturn}, we prove the results in \S\ref{ssec:introlatereturn}. \cref{sec:questions} contains open questions surrounding our results and a discussion thereof.

\subsection*{Acknowledgements}
We warmly thank Jean-Philippe Chassé, Jonny Evans, Rei Henigman, David Keren Yaar, Rémi Leclercq, Leonid Polterovich, Felix Schlenk, Sobhan Seyfaddini and Morgan Weiler for their interest, encouragement and their many stimulating comments on an earlier version which have greatly improved the present one. JB acknowledges the support of the following grants used during the completion of this project: Israel Science Foundation grant 1102/20, ERC Starting Grant 757585 and Swiss National Science Foundation Ambizione Grant PZ00P2-223460.
JK is supported by the National Research Foundation of Korea (NRF) grants funded by the Korean government (MIST, No. 2022R1F1A1074587) and through the G-LAMP program (MOE, RS-2024-00441954).

\section{Main theorem}
\label{sec:mainthm}

This section is dedicated to proving \cref{thm:main} and its refinement \cref{thm:main2}. 

\subsection{Symmetric probes}
\label{sec:symmprobes}

For the construction of equivalences of split tori, we make use of the standard toric structure $X_{\alpha} \rightarrow \square_{\alpha}$ for which split tori appear as toric fibres. This allows us to use symmetric probes first used for the purpose of constructing Hamiltonian isotopies between fibres. Let us recall the basics of this technique for the reader's convenience. For details, we refer to \cite{Bre23}. A \emph{symmetric probe} $\sigma \subset \Delta$ in the Delzant polytope arising as moment map image $\Delta \subset \R^n$ of a toric symplectic manifold is a segment intersecting $\pp \Delta$ integrally transversely in two points. Integral transversality means that $\sigma$ has a rational direction and that its primitive directional vector can be completed to a basis of the lattice by vectors parallel to the facet the probe intersects. This condition implies that one can perform symplectic $T^{n-1}$-reduction by a subgroup of the toric action and obtain a smooth copy of $S^2$ as quotient. Toric fibres corresponding to points lying on the symmetric probe map to circles in the reduced space, and Hamiltonian isotopies of these circles lift to Hamiltonian isotopies of the corresponding fibres. This yields \cite[Theorem A]{Bre23}. 

\begin{theorem}
\label{thm:symmprobes}
Let $X$ be a toric symplectic manifold with moment polytope $\Delta$ and $\sigma \subset \Delta$ a symmetric probe. If $p,q \in \sigma$ lie at equal distance to the boundary $\pp \Delta$, then $T(p) \cong T(q)$. 
\end{theorem}

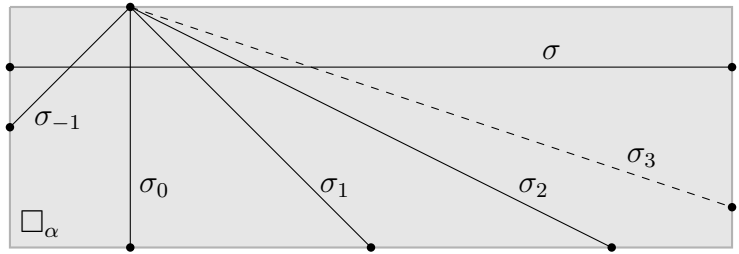
\begin{figure}
		\begin{tikzpicture}[scale=0.8]	
            \fill[black!10] (-6,2)--(6,2)--(6,-2)--(-6,-2)--(-6,2);
            \draw[thick,black!30] (-6,2)--(6,2)--(6,-2)--(-6,-2)--(-6,2);
            \draw[black] (-4,2)--(-6,0);
            \draw[black] (-4,2)--(-4,-2);
            \draw[black] (-4,2)--(-0,-2);
            \draw[black] (-4,2)--(4,-2);
            \draw[black, dashed] (-4,2)--(6,-1.33);
            \draw[black] (-6,1)--(6,1);
            
            \fill[thick] (-4,2)  circle[radius=2pt];
            \fill[thick] (-6,0)  circle[radius=2pt];
            \fill[thick] (-4,-2)  circle[radius=2pt];
            \fill[thick] (4,-2)  circle[radius=2pt];
            \fill[thick] (0,-2)  circle[radius=2pt];
            \fill[thick] (6,-1.33)  circle[radius=2pt];
            \fill[thick] (6,1)  circle[radius=2pt];
            \fill[thick] (-6,1)  circle[radius=2pt];

            \node at (-5.2,0.1){$\sigma_{-1}$};
            \node at (-3.6,-1){$\sigma_{0}$};
            \node at (-0.6,-1){$\sigma_{1}$};
            \node at (2.7,-1){$\sigma_{2}$};
            \node at (4.5,-0.5){$\sigma_{3}$};
            \node at (3,1.25){$\sigma$};
            \node at (-5.5,-1.6){$\square_{\alpha}$};
		\end{tikzpicture}
	\caption{Some examples and one non-example of symmetric probes in $\square_{\alpha}$. The segments $\sigma_k$ have slope $(k,-1)$ for $k \in \{-1,0,1,2,3\}$ and $\sigma$ has slope $(1,0)$. The dashed segment $\sigma_3$ is \emph{not} a symmetric probe, since its intersection with the the vertical edge is not integrally transverse.}
	\label{fig:2}
\end{figure}

Let us now discuss the case $X = X_{\alpha}$, where we have $\Delta = \square_{\alpha} = [-\alpha-1,\alpha+1] \times [-\alpha,\alpha] \subset \R^2$. The symmetric probes in $\square_{\alpha}$ are easy to understand, see \cref{fig:2}. Every line $\ell$ with slope~$(1,0)$ or~$(0,1)$ that intersects $\square_{\alpha}$ in its interior yields a symmetric probe. We call these probes \emph{horizontal} and \emph{vertical}, respectively. Every segment with slope~$(1,1)$ and~$(-1,1)$ which intersects~$\square_{\alpha}$ but none of its vertices yields a symmetric probe. The latter two types of probes will play a crucial role for us, as they generate Hamiltonian equivalences corresponding to billiards. Note that for small~$\alpha \in \R_{>0}$, there are other symmetric probes, namely those of slope~$(k,-1)$ for $k \in \Z$, which intersect only horizontal edges of~$\square_{\alpha}$. However, for~$\vert k \vert > 2$ these probes are redundant for the classification of toric fibres and determining their Hamiltonian monodromy group and will not be used in the paper.

\subsection{Billiards}
\label{ssec:billiards}

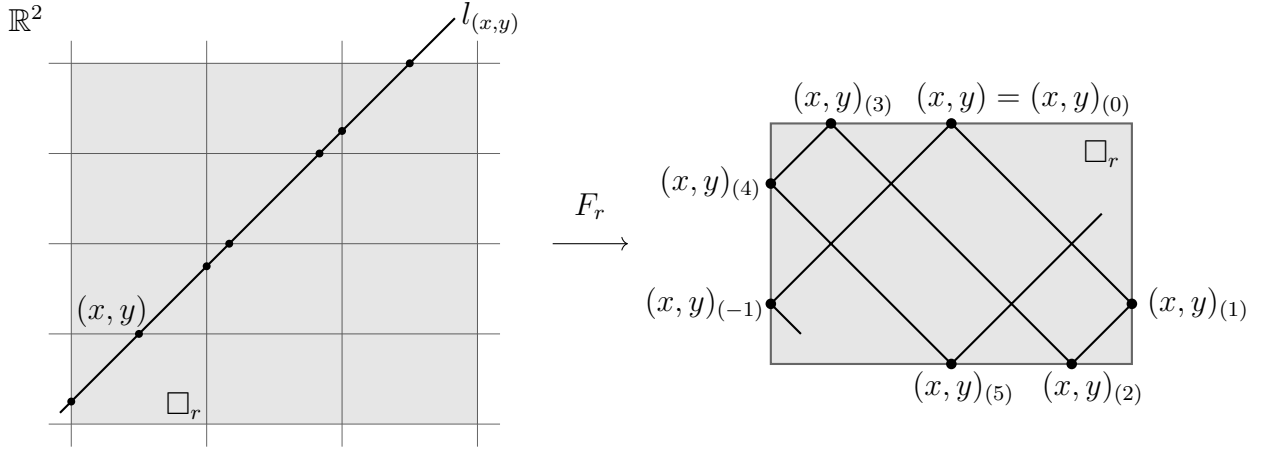
\begin{figure}
  \centering
  \begin{tikzpicture}
    \begin{scope}[scale=0.6,shift={(-8.5,0)}]
        \draw[very thin,black!60] (-3,-4.5)--(-3,4.5);
        \draw[very thin,black!60] (0,-4.5)--(0,4.5);
        \draw[very thin,black!60] (3,-4.5)--(3,4.5);
        \draw[very thin,black!60] (6,-4.5)--(6,4.5);
        \draw[very thin,black!60] (-3.5,-4)--(6.5,-4);
        \draw[very thin,black!60] (-3.5,-2)--(6.5,-2);
        \draw[very thin,black!60] (-3.5,0)--(6.5,0);
        \draw[very thin,black!60] (-3.5,2)--(6.5,2);
        \draw[very thin,black!60] (-3.5,4)--(6.5,4);
        \fill[black, opacity=0.1] (-3,-4) rectangle (6,4);
        \draw[thick,black] (-3.25,-3.75)--(5.5,5);
        \fill[thick] (-3,-3.5)  circle[radius=2.5pt];
        \fill[thick] (-1.5,-2)  circle[radius=2.5pt];
        \fill[thick] (0,-0.5)  circle[radius=2.5pt];
        \fill[thick] (0.5,0)  circle[radius=2.5pt];
        \fill[thick] (2.5,2)  circle[radius=2.5pt];
        \fill[thick] (3,2.5)  circle[radius=2.5pt];
        \fill[thick] (4.5,4)  circle[radius=2.5pt];
        \node at (6.3,5){$l_{(x,y)}$};
        \node at (-4,5){$\R^2$};
        \node at (-2.1,-1.5){$(x,y)$};
        \node at (-0.5,-3.6){$\square_r$};
    \end{scope}
    \draw[->] (-0.5,0) -- (0.5,0);
    \node at (0,0.5){$F_r$};
    \begin{scope}[scale=0.8,shift={(6,0)}]
        \fill[black, opacity=0.1] (-3,-2) rectangle (3,2);
        \draw[thick,black!60] (-3,-2) rectangle (3,2);
        \draw[thick,black] (-2.5,-1.5)--(-3,-1)--(0,2)--(3,-1)--(2,-2)--(-2,2)--(-3,1)--(0,-2)--(2.5,0.5);
        \fill[thick] (0,-2)  circle[radius=2.5pt];
        \fill[thick] (-3,-1)  circle[radius=2.5pt];
        \fill[thick] (0,2)  circle[radius=2.5pt];
        \fill[thick] (3,-1)  circle[radius=2.5pt];
        \fill[thick] (2,-2)  circle[radius=2.5pt];
        \fill[thick] (-2,2)  circle[radius=2.5pt];
        \fill[thick] (-3,1)  circle[radius=2.5pt];
        \node at (1.2,2.4){$(x,y)=(x,y)_{(0)}$};
        \node at (-4.1,-1){$(x,y)_{(-1)}$};
        \node at (4.1,-1){$(x,y)_{(1)}$};
        \node at (2.35,-2.4){$(x,y)_{(2)}$};
        \node at (-1.8,2.4){$(x,y)_{(3)}$};
         \node at (-4,1){$(x,y)_{(4)}$};
        \node at (0.2,-2.4){$(x,y)_{(5)}$};
        \node at (2.5,1.5){$\square_r$};
    \end{scope}
  \end{tikzpicture}
  \caption{Sketch of the folding map $F_r \colon \R^2 \rightarrow \square_r$. It maps the line $l_{(x,y)}$ of slope $(1,1)$ on the left-hand side to the good billiard trajectory of $(x,y)$ on the right hand side and intersection points with the grid to bouncing points.}
  \label{fig:3}
\end{figure}

Let $\alpha > 0$ and let $(x,y) \in \square_{\alpha}$. Recall that we have set $r(x,y) = \max\{\vert x \vert -1 , \vert y \vert\}$ and that the billiard table corresponding to $(x,y)$ is $\square_{r(x,y)}$ equipped with the equal angles reflection law. In the corner of the billiard table, the trajectory is reflected backwards. The \emph{good billiard trajectory} of $(x,y)$ was defined in \cref{def:goodbilliard}. Recall that a bouncing point of the good billiard trajectory of $(x,y) \in \square_{\alpha}$ is called \emph{admissible} if it does not coincide with a corner of the billiard table. In case $(x,y) \in \Sigma$, its good billiard trajectory is by definition stationary and has $(x,y)$ as single bouncing point, which we call admissible. We introduce some further notation. For any $k \in \Z$ let $(x,y)_{(k)} \in \pp \square_{r(x,y)}$ be the $k$-th \emph{admissible} bouncing point of the good billiard trajectory of $(x,y)$ with 
\begin{enumerate}
    \item $(x,y)_{(0)} = (x,y)$, \smallskip
    \item The clockwise branch of the billiard trajectory starting at $(x,y)_{(0)} = (x,y)$ corresponds to $k \geqslant 0$, the counterclockwise branch to negative $k \leqslant 0$.
\end{enumerate}
See also the right-hand side of \cref{fig:3}. For all $(x,y) \in \square_{\alpha}$, let
\begin{align*}
    B(x,y) 
    & = \{\text{admissible bouncing points of the good billiard trajectory of }(x,y)\} \\
    & = \{( x, y)_{(k)} \sth k \in \Z \} \subset \square_{r(x,y)}.
\end{align*}

The following proves the construction part of \cref{thm:main}.

\begin{lemma}
\label{prop:mainconstruction}
If some $(\pm x', \pm y') \in B(x,y)$ then $T(x,y) \cong T(x',y')$. Furthermore, this equivalence can be realized by a sequence of symmetric probes.
\end{lemma}

\proof All maps of the form $(x,y) \mapsto (\pm x, \pm y)$ can be realized by at most one horizontal and one vertical probe. If $(x,y) \in \Sigma$, there is nothing further to show, since $(x,y)_{(k)} = (x,y)$ for all $k$. Let $(x,y) \in \square_{\alpha} \setminus \Sigma$. We restrict our attention to the case of $k \geqslant 0$. The case $k \leqslant 0$ is similar. Let $n > 0$ be the smallest index such that the good billiard trajectory has hit a corner of $\pp \square_{r(x,y)}$ between $(x,y)$ and $(x,y)_{(n)}$. If such an index does not exist, set $n = + \infty$. Let $k < n$. We proceed by induction. The points $(x,y)_{(k)}$ and $(x,y)_{(k+1)}$ are not corner points of $\pp \square_{r(x,y)}$ and thus lie at equal distance to the boundary of a symmetric probe with direction $(1,1)$ or $(-1,1)$ and hence \cref{thm:symmprobes} implies the claim. For $k \geqslant n$, note that the billiard trajectory was reflected backwards at least once. Therefore, $(x,y)_{(k)} = (x,y)_{(k')}$ for some $k' < k$. This proves the claim. \proofend

To handle the combinatorics of billiards in the rectangle, we make use of the \emph{unfolding construction} for billiards, which we describe here. 
\begin{definition} For any $r > 0$, the \emph{folding map} is defined as
    \begin{equation}
        \label{eq:foldingmapdef}
        F_r \colon \R^2 \rightarrow  \square_r \subset \R^2, \quad
        (x,y) \mapsto ((-1)^m (x - 2m(r+1)) , (-1)^n (y - 2nr)),
    \end{equation}
    where $(m,n) \in \Z^2$ is chosen such that $F_r(x,y) \in \square_r$.
\end{definition}
Although the pair $(m,n) \in \Z^2$ is not unique whenever $x = (2l + 1)(1+r)$ or $y = (2l + 1)r$, the map is still well-defined. For later use, we also define, for $(m,n) \in \Z^2$, the following maps 
\begin{equation}
    A^r_{m,n} \colon \R^2 \rightarrow \R^2, \quad
    (x,y) \mapsto ((-1)^m x + 2m(r+1), (-1)^ny + 2nr).
\end{equation}
Note that 
\begin{equation}
    \label{eq:frarinverse}
    (F_r \circ A^r_{m,n})\vert_{\square_{\alpha}} = \id_{\square_{\alpha}},
\end{equation}
for all $m,n \in \Z$.
Under the folding map, the line $l_{(x,y)} = \{(x,y) +t(1,1) \sth t \in \R \} \subset \R^2$ is mapped to the good billiard trajectory of $(x,y)$ in $\square_r$. The set of admissible bouncing points $B(x,y)$ is given by
\begin{equation}
    \label{eq:goodbouncing}
    B(x,y) = F_r\left( l_{(x,y)} \cap F_r^{-1}(\pp\square_{r}) \setminus F_r^{-1}(\pm (1+r),\pm r) \right).
\end{equation}
Note that $F_r^{-1}(\pp\square_r)$ is the grid $\{ x = (2l_1 + 1)(1+r) \sth l_1 \in \Z \} \cup \{ y = (2l_2 + 1)r \sth l_2 \in \Z \}  \subset \R^2$ and $F_r^{-1}(\pm (1+r),\pm r)$ is the set of intersection points of the horizontal and vertical lines of that grid. These correspond to the corners of the billiard table, which is why we remove them. See \cref{fig:3} for a sketch of the folding map. Note that \eqref{eq:goodbouncing} holds, although $l_{(x,y)} \setminus F_r^{-1}(\pm (1+r),\pm r)$ may have several connected components (this happens exactly in case the good billiard trajectory of $(x,y)$ hits a corner of $\square_r$). Indeed, the images of different connected components of $l_{(x,y)} \setminus F_r^{-1}(\pm (1+r),\pm r)$ under $F_r$ coincide.

For $\alpha>0$, let

\begin{equation}
    \label{eq:Qdefset}
    Q = \{ (x,y) \in \R^2 \sth y >0 , y > \vert x \vert - 1 \},
\end{equation}
see also \cref{fig:1}. It is enough to classify split tori in $(Q \cup \Sigma) \cap \square_{\alpha} \subset \square_{\alpha}$. Indeed, any split torus in the complement of this set is equivalent by a single symmetric probe to a split torus in $Q$. We thus restrict our attention to $Q$ to simplify the statement of the following lemma. Without this restriction, there would be a distinction of cases in \eqref{eq:billiard0} according to which expression realizes the maximum $r(x,y) = \max\{\vert x \vert -1, \vert y \vert \}$.

\begin{lemma}
\label{lem:billiard}
Let $\alpha > 0$ and $(x,y),(x',y) \in Q$. If there are $k_1,k_2 \in \Z$ and $\delta \in \{0,1\}$ with
\begin{equation}	
	\label{eq:billiard0}
	x =  (-1)^{\delta} x' + 2k_1 + 2k_2y,
\end{equation}
then there are $\delta_1,\delta_2 \in \{0,1\}$ such that $((-1)^{\delta_1} x' , (-1)^{\delta_2} y) \in B(x,y)$. 
\end{lemma}

Note that if $(x,y),(x',y') \in \Sigma$, then, for some choice of signs, we have $(\pm x' , \pm y') \in B(x,y) = \{(x,y)\}$.\smallskip

\proofof{Lemma~\ref{lem:billiard}}
Since $(x,y) \in Q$, the billiard table is $\square_r$ for $r=y$. We use the folding construction to prove the statement. In particular, recall from \eqref{eq:goodbouncing} that points in $B(x,y)$ lift to intersection points of the grid $F_r^{-1}(\pp \square_r)$ with $l_{(x,y)}$ under the unfolding. Since $(x,y),(x',y)$ are not on the corner of $\square_r$, their lifts are not in $F_r^{-1}(\pm(1+r),\pm r)$. Furthermore, since $(x,y),(x',y) \in Q$ we can restrict our attention to intersection points of $l_{(x,y)}$ with horizontal lines of the grid. These are given by the set 
\begin{equation}
    \label{eq:set1}
    \{ (x,y) + 2kr(1,1) \sth k \in \Z \}.
\end{equation}
By \cref{prop:mainconstruction}, the image $F_r(x'',y'')$ under the folding map of every $(x'',y'') \in \R^2$ which is in the set \eqref{eq:set1} but not in $F_r^{-1}(\pm(1+r),\pm r)$ satisfies $T(x,y) \cong T(F_r(x'',y''))$ by applying a billiard trajectory of symmetric probes. On the other hand, the set of points in $\R^2$ which map to $(\pm x', \pm y)$ under the folding map is equal to
\begin{equation}
    \label{eq:set2}
    F_r^{-1}(\pm x', \pm y) 
    =
    \left\{ A_{m,n}^r((-1)^{\delta_1}x', (-1)^{\delta_2}y) \sth (m,n) \in \Z^2, \delta_1,\delta_2 \in \{0,1\} \right\},
\end{equation}
where we have used \eqref{eq:frarinverse}. Therefore, to prove the statement it is sufficient to prove that the intersection of the sets \eqref{eq:set1} and \eqref{eq:set2} is non-empty. A computation shows that setting $m = k_1$, $k = n = k_1 - k_2$ and taking $\delta_1,\delta_2$ such that $\delta_1 \equiv \delta + m \mod 2$ and $\delta_2 \equiv n \mod 2$ yields a point lying in both sets.
\proofend

\subsection{Proof of main theorem}

\cref{thm:main} is a consequence of the following. 

\begin{theorem}
\label{thm:main2}
Let $\alpha > 0$ and let $(x,y),(x',y') \in Q \cup \Sigma$, then the following are equivalent:
\begin{enumerate}
    \item One of the points $(\pm x', \pm y')$ is in the set $B(x,y)$ of admissible bouncing points of the good billiard trajectory of $(x,y)$;
    \item The split tori $T(x,y)$ and $T(x',y')$ are Hamiltonian isotopic;
    \item If $(x,y) \in \Sigma$, then $(x',y') \in \Sigma$ and $(x,y)$ is equal to some $(\pm x', \pm y')$. If $(x,y), (x',y') \in Q$, then 
    \begin{equation}
        y' = y, \quad x = \pm x' + 2k_1 + 2k_2y, \quad \text{for some } k_1,k_2 \in \Z.
    \end{equation}
\end{enumerate}
\end{theorem}

\cref{prop:mainconstruction} shows (1) $\Rightarrow$ (2) and \cref{lem:billiard} shows (3) $\Rightarrow$ (1). To complete, we need to show (2) $\Rightarrow$ (3), which is \cref{lem:mainobstruction}. To prepare its proof, let us discuss some of its ingredients.

We use obstructions from \cite{Bre23}, which apply to toric fibres in general compact toric manifolds. Let $X$ be a compact toric symplectic manifold with moment polytope $\Delta \subset \R^n$ and denote by $\ell_i(p)$ for $1 \leqslant i \leqslant N = \#\{\text{facets of } \Delta \}$ the integral affine distance of a point $p \in \Int \Delta$ to the $i$-th facet of $\Delta$. The following is \cite[Theorem B]{Bre23}.

\begin{theorem}
\label{thm:Chekinvts}
Let $X$ be a compact toric symplectic manifold with moment polytope $\Delta$. If $T(p) \cong T(q)$ for $p,q \in \Int \Delta$, then 
\begin{enumerate}
    \item $d(p) = d(q)$, where $d(p) = \min_i \{\ell_i(p)\}$,
    \item $\#_d(p) = \#_d(q)$, where $\#_d(p) = \{ i \sth \ell_i(p) = d(p) \}$,
    \item $\Gamma(p) = \Gamma(q)$, where $\Gamma(p) = \Z\langle \ell_1(p) - d(p), \ldots, \ell_N(p) - d(p) \rangle \subset \R$.
\end{enumerate}
\end{theorem}

Here, we denote by $\Z \langle a_1,\ldots,a_n \rangle \subset \R$ the subgroup of linear combinations with integer coefficients of $a_1,\ldots,a_n \in \R$. These invariants are derived, via Delzant's construction, from the invariants for product tori from \cite{Che96}, which is why we call them \emph{Chekanov invariants}. We refer to Sections 1.2 and 4.1 in \cite{Bre23} for a thorough discussion. Though these invariants are not complete, we will use the first and second Chekanov invariant in the proof of \cref{lem:mainobstruction}. For $\Delta = \square_{\alpha}$, the functions $\ell_i$ for $i \in \{1,2,3,4\}$ are given by
\begin{equation}
\ell_1(x,y) = \alpha - y, \quad 
\ell_2(x,y) = 1 + \alpha - x,\quad
\ell_3(x,y) = \alpha + y, \quad
\ell_4(x,y) = 1 + \alpha + x.
\end{equation}
We also use \cite[Theorem 4.7]{Bre23}, which yields constraints on the map $\phi_* \colon H_2(X,T(p)) \rightarrow H_2(X,T(q))$ induced on relative homology by $\phi \in \Ham(X)$ mapping a toric fibre $T(p)$ to another toric fibre $T(q)$. We refer to \cite[Section 4.2]{Bre23} for details. 

\begin{lemma}
\label{lem:mainobstruction}
Let $\alpha > 0$. If $(x,y) \in \Sigma \subset \square_{\alpha}$, and $T(x,y) \cong T(x',y')$ then $(x,y) = (\pm x',\pm y')$. If $(x,y),(x',y') \in Q$ and $T(x,y)\cong T(x',y')$ then 
\begin{equation}
    y = y' \text{ and } x = \pm x' + 2k_1 + 2k_2 y \text{ for some }k_1,k_2 \in \Z.
\end{equation}
\end{lemma}

\proof Note that $\Sigma \subset \square_{\alpha}$ is the subset where the second Chekanov invariant is strictly greater than one. In particular, split tori in $\Sigma$ can only be equivalent to split tori in $\Sigma$, which allows us to treat this case separately. The only points having $\#_d(x,y) = 3$ are $(0,1)$ and $(0,-1)$, which proves the claim for these. Let $(x,y)$ be any other point of $\Sigma$, then $d(x,y) = \alpha - \vert y \vert$ and $\#_d(x,y) = 2$. In case $\vert y \vert > 0$, this allows us to conclude, since the first Chekanov invariant yields $\vert y \vert = \vert y' \vert$ and thus $\vert x \vert = \vert x' \vert$. In case $y = 0$ (this corresponds to the segment $[-1,1] \times \{0\}$), we cannot yet conclude and it will be treated similarly to the case of $(x,y) \in Q$, see the end of the proof. Alternatively, it can be settled by appealing to \cite{FOOO13}, as these product tori are distinguished by their displacement energy.

Now let $(x,y),(x',y') \in Q$. Then $y=y'$ follows from the first Chekanov invariant, $d(x,y)=\alpha - y$. The second Chekanov invariant is constant equal to one on $Q$ and therefore of no further use. The third Chekanov invariant is not complete and we use \cite[Theorem 4.7]{Bre23} instead. Assume there is~$\phi \in \Ham(X_{\alpha})$ mapping~$L = T(x,y)$ to~$L' = T(x',y)$. Then we have the following commutative diagram
\begin{equation}
\label{eq:homologycd}
	\begin{tikzcd}
		0 \arrow[r] 
		& H_2(X_{\alpha}) \arrow[r] \arrow[d, "\id"]
		& H_2(X_{\alpha},L) \arrow[r,"\partial"] \arrow[d, "\phi_*"]
		& H_1(L) \arrow[r] \arrow[d, "(\phi\vert_L)_*"]
		& 0 \\		
		0 \arrow[r] 
		& H_2(X_{\alpha}) \arrow[r]
		& H_2(X_{\alpha},L') \arrow[r,"\partial'"]
		& H_1(L') \arrow[r]
		& 0.
	\end{tikzcd}
\end{equation}
As the inclusion~$L \subset X_{\alpha}$ yields the zero map in homology, this diagram can be easily extracted from the long exact sequences on relative homology. On $H_2(X_{\alpha})$ the Hamiltonian diffeomorphism~$\phi$ induces the identity since~$\phi$ is isotopic to the identity. We use~\cite[Theorem 4.7]{Bre23} to determine which isomorphisms~$H_2(X_{\alpha},L) \rightarrow H_2(X_{\alpha},L')$ can be induced by the Hamiltonian diffeomorphism~$\phi$. Let~$D_1,D_2,D_3,D_4 \in H_2(X_{\alpha},L)$ be the classes defined by the obvious disks bounding the circles in the product $S^2 \times S^2$, labelled such that $D_1$ and $D_3$ contain the north and south poles of the smaller sphere, respectively, and $D_2$ and $D_4$ the north and south poles of the larger sphere, respectively. This notation is coherent with the one used in~\cite[Proposition 2.9]{Bre23} for the general toric case. All of the classes~$D_i$ are represented by topological disks of Maslov index~$2$ and have symplectic area equal to the integral affine distance of~$(x,y) \in \Int \square_{\alpha}$ to the corresponding facet of~$\square_{\alpha}$, i.e.\ $\int_{D_i} \omega_{\alpha} = \ell_i(x,y)$. In our case, we have, 
\begin{equation}
	\int_{D_1} \omega_{\alpha} = \alpha - y, \quad
	\int_{D_2} \omega_{\alpha} = 1 + \alpha - x, \quad
	\int_{D_3} \omega_{\alpha} = \alpha + y, \quad
	\int_{D_4} \omega_{\alpha} = 1 + \alpha + x.
\end{equation}
Similarly, we denote by~$D_i'$ the corresponding classes in~$H_2(X_{\alpha},L')$. Recall from~\cite[Definition 4.6]{Bre23} that the set of distinguished classes~$\cd(x,y) \subset H_2(X_{\alpha},L)$ is formed by the~$D_i$ for which the integral affine distance to the boundary~$\pp \square_{\alpha}$ is minimal. Since both base points are in $Q$, we have~$\cd(x,y) = \{D_1\}$ and~$\cd(x',y) = \{D_1'\}$, meaning that~$\phi_* D_1 = D_1'$ by~\cite[Theorem 4.7]{Bre23}. Now set
\begin{equation}
    \label{eq:phistarrr}
	\phi_* D_2
	=
	a_1 D_1' + a_2 D_2' + a_3 D_3' + a_4 D_4', \quad
	a_i \in \Z. 
\end{equation}
The integers~$a_i$ determine~$\phi_*$. Indeed, note that~$D_1 + D_3, D_2 + D_4 \in H_2(X_{\alpha})$ and therefore~$\phi_*(D_1 + D_3) = D_1' + D_3'$ and~$\phi_*(D_2 + D_4) = D_2' + D_4'$. We deduce that~$\phi_*D_3 = D_3'$ and~$\phi_*D_4 = -a_1D_1' + (1-a_2) D_2' - a_3 D_3' + (1-a_4)D_4'$. Since~$\phi_*$ is an isomorphism, we obtain
\begin{equation}
	\label{eq:detphi}
	\det \phi_*
	= 
	a_2 - a_4 
	= 
	\pm 1.
\end{equation}
Let us first treat the case~$\det \phi_* = +1$. Then~$a_4 = a_2 - 1$. Since~$\phi_*$ preserves the Maslov class, \eqref{eq:phistarrr} yields~$a_3 = 2 - a_1  - 2a_2$. Since it preserves the area class, we find
\begin{eqnarray*}
	1 + \alpha - x 	
	&=& \int_{D_2} \omega_{\alpha} \\
	&=& \int_{\phi_* D_2} \omega_{\alpha} \\
	&=& a_1(\alpha - y) + a_2(1 + \alpha - x') + (2-a_1-2a_2)(\alpha + y) + (a_2 - 1)(1 + \alpha + x') \\
	&=& 2a_2 + 2(1 - a_1 - a_2)y - 1 + \alpha - x',
\end{eqnarray*}
from which we deduce
\begin{equation}
	\label{eq:main+1}
	x  
	=
	x' + 2(1 - a_2) + 2(a_1 + a_2 -1)y
	=
    x' + 2k_1 + 2k_2y,
\end{equation}
for $k_1 = 1 - a_2$ and $k_2 = a_1 + a_2 - 1$. In case $\det \phi_* = -1$, an analoguous computation shows 
\begin{equation}
	\label{eq:main-1}
	x 
	=
	-x' -2a_2 + 2(a_1 + a_2)y
	=
    -x' + 2k_1 + 2k_2y,
\end{equation}
for $k_1 = -a_2$ and $k_2 = a_1 + a_2$. This proves the statement for $(x,y),(x',y') \in Q$. 

It remains to show the statement for $(x,y),(x',y') \in [-1,1] \times \{0\}$, in which case $\cd(x,0) = \{D_1,D_3\}$, and $\cd(x',0) = \{D_1',D_3'\}$. This means that either~$\phi_* D_1 = D_1'$ or~$\phi_* D_1 = D_3'$. In the former case, the same computation as in the case of $Q$ can be carried out. In the case~$\phi_* D_1 = D_3'$, a similar argument applies, up to a change in sign in~\eqref{eq:detphi}. This means that in both cases we also obtain equations~\eqref{eq:main+1} and~\eqref{eq:main-1}. Since $y=0$ in this case, this completes the proof.\proofend

\section{Hamiltonian monodromy}

\label{sec:monodromy}

The goal of this section is to discuss the Hamiltonian monodromy group of split tori in $X_{\alpha}$ and prove our main result in that direction, \cref{thm:monodromy}.

\subsection{Preparations}
\label{ssec:recollections_monodromy}
The Hamiltonian monodromy group is defined in \eqref{eq:hammonodromygroup}. Note that it is a symplectic invariant up to conjugacy in the sense that for all $\phi \in \Symp(X,\omega)$, we have
\begin{equation}
	\label{eq:monodromyinv}
	\ch_{\phi(L)} = (\phi\vert_L)_* \circ \ch_L \circ (\phi\vert_L)_*^{-1}.
\end{equation}
Thus it is enough to determine $\ch_{L}$ for a representative of a Hamiltonian orbit in the space of split tori. 

Let $\sigma \subset \square_{\alpha}$ be a symmetric probe with primitive directional vector $v \in \Z^2$ intersecting facets defined by normal vectors $\xi,\xi' \in \Z^2$. Here we choose $\xi,\xi' \in \Z^2$ to be primitive and pointing into the interior of $\square_{\alpha}$ with the convention that $v$ points from the facet defined by $\xi$ to the one defined by $\xi'$. Let $L,L' \subset X_{\alpha}$ split tori lying on $\sigma$ at equal distance to the boundary. By \cref{thm:symmprobes}, the symmetric probe induces a Hamiltonian diffeomorphism $\phi \in \Ham (X_{\alpha})$ with $\phi(L) = L'$. Via the toric group action, there is a natural way to identify $H_1(L) \cong H_1(L')$, which we will use tacitly from now on, and to interpret $v \in H^1(L)$ and $\xi,\xi' \in H_1(L)$. With this interpretation, the Hamiltonian diffeomorphism $\phi$ induces the following map on the first homology, 
\begin{equation}
	\label{eq:probehomology}
	   (\phi\vert_L)_* \colon H_1(L) \rightarrow H_1(L'), \quad
	   a \mapsto a + \langle v , a \rangle (\xi' - \xi).
\end{equation}
For details, we refer to \cite[Section 3]{Bre23}, in particular to Theorem 3.2 and Figure 4 therein. In the remainder of this section, we express automorphisms of the first homology of split tori by two-by-two matrices expressed in the basis on $H_1(L)$ induced by the Hamiltonian $T^2 = \R^2/\Z^2$-action on $X_{\alpha}$ via the moment map $X_{\alpha} \rightarrow \square_{\alpha}$. Up to choosing adequate orientations, this is the same basis as the natural one obtained from viewing a split torus $T(x,y) = S^1_x \times S^1_y \subset S^2 \times S^2$ as product of circles and taking as basis the elements $S^1_x \times \{*\}$ and $\{*\} \times S^1_y$. Recall that, in this paper, we only use symmetric probes with directional vectors $(1,0),(0,1),(1,1),(-1,1)$ and, by \eqref{eq:probehomology}, these yield the six matrices
\begin{equation}
    \label{eq:matrices}
    \begin{pmatrix}
        - 1 & 0 \\
        0 & 1
    \end{pmatrix},
    \begin{pmatrix}
        1 & 0 \\
        0 & -1
    \end{pmatrix},
    \begin{pmatrix}
		1 & 0 \\
		\pm 2 & -1
	\end{pmatrix}, 
    \pm
    \begin{pmatrix}
		0 & 1 \\
		  1 & 0
	\end{pmatrix}
\end{equation}
in the above basis. The first two matrices come from the horizontal and vertical probe, respectively. The following pair of matrices comes up in the case of symmetric probes $(- 1,1)$ and $(1,1)$ having endpoints on two horizontal edges. The latter pair of matrices comes up in case the same symmetric probes have endpoints on one horizontal and one vertical edge. 

\subsection{Obstructions}

We prove the following upgrade of \cref{lem:mainobstruction}.

\begin{lemma}
\label{lem:mainobstructionmonodromy}
Let $L=T(x,y)$ and $L'=T(x',y)$ for $(x,y),(x',y) \in Q$. If there is $\phi \in \Ham (X_{\alpha})$ with $\phi(L)=L'$ then the map $(\phi\vert_L)_* \colon H_1(L) \rightarrow H_1(L')$ has the form
\begin{equation}
    \label{eq:monodromyQ}
    (\phi\vert_L)_*
    =
    \begin{pmatrix}
        (-1)^{\delta} & 0 \\
        2k_2 & 1 
    \end{pmatrix}
\end{equation}
for some $\delta \in \{0,1\}$ and $k_1, k_2 \in \Z$ satisfying
\begin{equation}
    \label{eq:monodromyQ2}
    x = (-1)^{\delta}x' + 2k_1 + 2k_2y.
\end{equation}
\end{lemma}

\proof 
This follows from the proof of \cref{lem:mainobstruction}. Indeed, consider the map $\phi_* \colon H_2(X_{\alpha},L) \rightarrow H_2(X_{\alpha},L')$ induced on relative second homology and let $\phi_* D_2 = \sum_i a_i D_i'$. Furthermore let $\delta \in \{0,1\}$ such that $\det \phi_* = (-1)^{\delta}$. For $\delta = 0$, we obtain $x = x' + 2(1-a_2) + 2(a_1 + a_2 - 1)y$ as in \eqref{eq:main+1}. Setting $k_1 = 1-a_2$ and $k_2 = a_1 + a_2 -1$ and using $\pp \circ \phi_* = (\phi\vert_L)_* \circ \pp$ as in \eqref{eq:homologycd}, we find 
\begin{equation}
     \begin{pmatrix}
        1 & 0 \\
        2k_2 & 1 
    \end{pmatrix}.
\end{equation}
The case $\delta = 1$ is treated similarly using \eqref{eq:main-1} and yields 
\begin{equation}
     \begin{pmatrix}
        -1 & 0 \\
        2k_2 & 1 
    \end{pmatrix}.
\end{equation}
\proofend

Note that for $\delta=1$, \eqref{eq:monodromyQ} yields an involution, whereas for $\delta=0$, it yields an element of infinite order.

\subsection{Hamiltonian monodromy and billiards}

Returning to the folding construction of billiards, we show the following refinement of \cref{lem:billiard}. Contrary to \cref{prop:mainconstruction} and \cref{lem:billiard}, we cannot circumvent the fact that some trajectories hit corners of the billiard table, see \cref{ex:badtorus}.

\begin{lemma}
\label{lem:billiardmonodromy}
Let $0<y<\alpha$ and $(x,y),(x',y) \in Q$ and let $k_1,k_2 \in \Z$ and $\delta \in \{0,1\}$ with
\begin{equation}	
	x =  (-1)^{\delta} x' + 2k_1 + 2k_2y.
\end{equation}
Assume furthermore that the corresponding billiard trajectory between $(x,y)$ and $(\pm x', \pm y)$ does not hit any corners of the billiard table. Then there is a Hamiltonian diffeomorphism $\phi \in \Ham (X_{\alpha})$ obtained by concatenation of symmetric probes mapping $L = T(x,y)$ to $L' = T(x',y)$ such that the induced map $(\phi\vert_L)_* \colon H_1(L) \rightarrow H_1(L')$ satisfies
    \begin{equation}
	\label{eq:billiard}
	(\phi \vert _L)_* 
	= 
	\begin{pmatrix}
		(-1)^{\delta} & 0 \\
		2k_2 & 1
	\end{pmatrix}.
\end{equation}
\end{lemma}

By \emph{the corresponding billiard trajectory between $(x,y)$ and $(\pm x', \pm y)$} we mean the one from \cref{lem:billiard}. From the same lemma, it follows that there is such a $\phi \in \Ham (X_{\alpha})$. What is new here is the claim about its monodromy \eqref{eq:billiard}. \smallskip

\proofof{\cref{lem:billiardmonodromy}} Since $(x,y),(x',y) \in Q$, we have $r=r(x,y) = r(x',y) = y$. For the proof, we assume $k_1 - k_2 \geqslant 0$. The case $k_1 - k_2 < 0$ is treated similarly.
Recall from the proof of \cref{lem:billiard} that
\begin{equation}
    (x,y) + 2ky(1,1) = A_{m,n}^y((-1)^{\delta_1}x',(-1)^{\delta_2}y),
\end{equation}
for $m = k_1$, $n = k = k_1 - k_2$ and $\delta_1,\delta_2$ such that $\delta_1 \equiv \delta + m \mod 2$ and $\delta_2 \equiv n \mod 2$. By assumption, the billiard trajectory does not hit any corners of the table $\square_r$. Let $\chi \in \Ham(X_{\alpha})$ be a Hamiltonian diffeomorphism generated by the symmetric probes corresponding to the billiard trajectory. Then $\chi(L) = T((-1)^{\delta_1}x',(-1)^{\delta_2}y)$. Let us compute $(\chi\vert_L)_*$. Since $n=k$, this billiard has $((-1)^{\delta_1}x',(-1)^{\delta_2}y)$ as bouncing point after hitting \emph{horizontal} edges of $\pp\square_r$ $k$ times. Let $i \in \{1, \ldots, k\}$ and let $\psi_i$ be the map induced on $H_1(\cdot)$ of the split tori corresponding to the bouncing points between the $(i-1)$-st and the $i$-th bouncing point on a horizontal edge. Then 
\begin{equation}
    \label{eq:chiexpansion}
    (\chi\vert_L)_* = \psi_k \circ \psi_{k-1} \circ \ldots \circ \psi_1. 
\end{equation}
For every $i$ there are two possibilities. Either there is a single bouncing point (on a vertical edge of the billiard table) in between the $(i-1)$-st and the $i$-th horizontal bouncing point, or there is not. In the former case, its billiard trajectory is a concatenation of two probes and, using \eqref{eq:matrices} it acts by
\begin{equation}
	\label{eq:phitype1}
	\psi_i
	=
	-\id
	=
	\begin{pmatrix}
		-1 & 0 \\
		0 & -1
	\end{pmatrix}. 
\end{equation}
on homology. We say that $\psi_i$ is of \emph{first type}. In the latter case, its billiard trajectory is realized by a single symmetric probe acting on homology by 
\begin{equation}
	\label{eq:phitype2}
	\psi_i
	=
	\begin{pmatrix}
		1 & 0 \\
		\pm 2 & -1
	\end{pmatrix},
\end{equation}
where the sign depends on whether the induced symmetric probe has directional vector~$(1,-1)$ or directional vector~$(1,1)$. See again \eqref{eq:probehomology} and \eqref{eq:matrices}. We say that $\psi_i$ is of \emph{second type}. Since $m = k_1$, the billiard trajectory hits $k_1$ vertical edges and thus the matrix of first type appears $k_1$ times in the expression \eqref{eq:chiexpansion}. The $\pm$-signs in the matrices of second type appear in alternating order starting with a $+$, since the billiard trajectory starts in the clockwise direction. It is here where we use $k_1 - k_2 = k \geqslant 0$, since otherwise the the billiard trajectory starts in the counterclockwise direction and the product of matrices starts with a $-$ sign. The matrices of second type appear $k-k_1 = -k_2$ times\footnote{The notation is inconvenient here; actually $-k_2$ is nonnegative. Indeed, we are in the case where $k_1-k_2 >0$ and thus $0 \leqslant m = k_1 \leqslant n = k_1 - k_2$, whence $k_2 \leqslant 0$. Here we have used $m \leqslant n$, which can be easily deduced from the folding picture of billiards.}. Since the matrix of the first type is $-\id$, it commutes with every other matrix and we can write
\begin{equation}	
	\label{eq:chiexpression}
	(\chi \vert _L)_*
	= 
	(-\id)^{k_1}
	\begin{pmatrix}
		1 & 0 \\
		(-1)^{k_2}2k_2 & (-1)^{k_2}
	\end{pmatrix}
	= 
	\begin{pmatrix}
		(-1)^{k_1} & 0 \\
		(-1)^{k_1 + k_2}2k_2 & (-1)^{k_1 + k_2}
	\end{pmatrix},
\end{equation}
where the first factor in the product comes from applying the matrix of first type $k_1$ times and the second factor from applying matrices of the second type $-k_2$ times. To compute the matrix $(\phi\vert_L)_*$, we post-compose $\chi$ with swaps of the sphere factors to compensate for the signs in $\chi(L) = T((-1)^{\delta_1}x',(-1)^{\delta_2}y)$. This yields
\begin{equation*}	
	\label{eq:phiexpression}
	(\phi \vert _L)_*
	= 
	\begin{pmatrix}
		(-1)^{\delta_1} & 0 \\
		  0 & (-1)^{\delta_2}
	\end{pmatrix}
	\begin{pmatrix}
		(-1)^{k_1} & 0 \\
		(-1)^{k_1 + k_2}2k_2 & (-1)^{k_1 + k_2}
	\end{pmatrix}
    =
    \begin{pmatrix}
		(-1)^{k_1 + \delta_1} & 0 \\
		(-1)^{k_1 + k_2 + \delta_2}2k_2 & (-1)^{k_1 + k_2 + \delta_2}
	\end{pmatrix}.
\end{equation*}
By the above choices, we find $k_1 + \delta_1 = m + \delta_1 \equiv \delta \mod 2$ and $\delta_2 \equiv n = k_1 - k_2 \equiv k_1 + k_2 \mod 2$ and thus \eqref{eq:billiard} follows.

The case $k_1 - k_2 < 0$ is proved similarly, since \eqref{eq:chiexpression} holds in that case, too. \proofend

The condition that the billiard trajectory between $(x,y)$ and $(\pm x', \pm y)$ does not hit a corner of the billiard table $\square_r$ is difficult to get a handle on. For simplicity, we restrict our attention to the set of points whose billiard trajectory \emph{never} hits a corner point. 

\begin{definition}
By $S_r \subset \pp \square_r$ we denote the set of points whose billiard trajectory eventually hits a corner.
\end{definition}

In the notation of \cref{ssec:billiards},
\begin{equation}
    S_r = \left\{ (x,y) \in \pp\square_r \sth l_{(x,y)} \cap F_r^{-1}(\pm (1+r) , \pm r) \neq \varnothing \right\}.
\end{equation}

Together with \cref{lem:mainobstructionmonodromy}, \cref{lem:billiardmonodromy} has the following immediate consequence. 

\begin{corollary}
\label{cor:monodromyiff}
Let $0<y<\alpha$ and $L = T(x,y), L' = T(x',y)$ for $(x,y),(x',y) \in Q \setminus S_y$. Then there is $\phi \in \Ham (X_{\alpha})$ with $\phi(L) = L'$ with $(\phi\vert_L)_* = \Phi$ for an isomorphism $\Phi \colon H_1(L) \rightarrow H_1(L')$ if and only if 
 \begin{equation}
	\label{eq:billiardmon}
	\Phi
	= 
	\begin{pmatrix}
		(-1)^{\delta} & 0 \\
		2k_2 & 1
	\end{pmatrix}
\end{equation}
for some $\delta \in \{0,1\}$ and $k_1,k_2 \in \Z$ satisfying $x = (-1)^{\delta}x' + 2k_1 + 2k_2y$.
\end{corollary}

The intersection of $S_r$ with $Q$ (where $r=y$) has a simple description which will be used below, 
\begin{equation}    
    \label{eq:badx}
    S_y \cap Q = \{(x,y) \in Q \sth x = 2l_1 + 1 + (2l_2 + 1)y \text{ for some } l_1,l_2 \in \Z \}.
\end{equation}
Indeed, the billiard trajectory of $(x,y) \in \pp \square_y \cap Q$ hits one of the corner points $(\pm (1+y), \pm y)$ if and only if there are $k_1,k_2 \in \Z$ for which $x = \pm(1+y) + 2k_1 + 2k_2y$ which is the case if and only if $(x,y)$ is in the set on the right hand side of \eqref{eq:badx}. This follows from the same methods as used in the proof of \cref{lem:billiard} by setting $x' = 1+y$.

\subsection{Monodromy theorem} \label{ssec:monodromythm} The goal of this section is to determine the Hamiltonian monodromy group of toric fibres in~$X_{\alpha}$. In light of~\eqref{eq:monodromyinv}, it is sufficient to compute the Hamiltonian monodromy group of~$T(x,y)$ for~$(x,y) \in \Sigma \cup Q$ as defined in \eqref{eq:Sigma} and \eqref{eq:Qdefset}. Subdivide $\Sigma$ further as $\Sigma = \{0,0\} \sqcup \Sigma_1 \sqcup \{(\pm 1, 0)\} \sqcup \Sigma_2$, where $\Sigma_1 = \Sigma \cap \{0 < \vert x \vert < 1\}$ and $\Sigma_2 = \Sigma \cap \{\vert x \vert > 1\}$. We recall that the obstruction from \cref{lem:mainobstructionmonodromy} only holds for $(x,y) \in Q$, meaning there is no contradiction with the fact that the Hamiltonian monodromy group for some $(x,y) \in \Sigma$ is larger. 

\begin{theorem}
\label{thm:monodromy}
For~$(x,y) \in \Sigma$ the Hamiltonian monodromy group of $L = T(x,y)$ is given by
\begin{equation}
	\label{eq:mon1.5}
	\ch_{L} 
	=
	\left\langle
	\begin{pmatrix}
		1 & 0 \\
		2 &  1
	\end{pmatrix},
	\begin{pmatrix}
		1 & 0 \\
		0 & - 1
	\end{pmatrix},
	\begin{pmatrix}
		-1 & 0 \\
		0 & 1
	\end{pmatrix}
	\right\rangle
	= 
	\left\{ 
	\left.
	\begin{pmatrix}
		\pm 1 & 0 \\
		2k & \pm 1
	\end{pmatrix}
	\right\vert
	k \in \Z
	\right\} 
	\cong \Z_2 \ltimes \Z \rtimes \Z_2,
\end{equation}
if $(x,y) = (0,0)$, by
\begin{equation}
	\label{eq:mon1}
	\ch_L 
	=
	\left\langle
	\begin{pmatrix}
		1 & 0 \\
		2 &  1
	\end{pmatrix},
	\begin{pmatrix}
		1 & 0 \\
		0 & - 1
	\end{pmatrix}
	\right\rangle
	= 
	\left\{ 
	\left.
	\begin{pmatrix}
		1 & 0 \\
		2k & \pm 1
	\end{pmatrix}
	\right\vert
	k \in \Z
	\right\} 
	\cong \Z_2 \ltimes \Z,
\end{equation}
if~$(x,y) \in \Sigma_1$, by 
\begin{equation}
	\label{eq:mon3}
	\ch_{L} 
	=
	\left\langle
	\begin{pmatrix}
		1 & 0 \\
		0 & -1
	\end{pmatrix}
	\right\rangle
	= 
	\left\{ 
	\begin{pmatrix}
		1 & 0 \\
		0 & \pm 1
	\end{pmatrix}
	\right\} \cong \Z_2.
\end{equation}
if $(x,y) = (\pm 1, 0)$, and by
\begin{equation}
	\label{eq:mon2}
	\ch_{L}
	=
	\left\langle
	\begin{pmatrix}
		0 & 1 \\
		1 & 0
	\end{pmatrix}
	\right\rangle
	=
	\left\{ 
	\begin{pmatrix}
		1 & 0 \\
		0 & 1
	\end{pmatrix},
	\begin{pmatrix}
		0 & 1 \\
		1 & 0
	\end{pmatrix}
	\right\} \cong \Z_2,
\end{equation} 
if~$(x,y) \in \Sigma_2$.

For the remainder of the statement, let $(x,y) \in Q$. First, we treat the case where $y$ is irrational. If furthermore $x \notin \Z\langle 1,y \rangle$, then the Hamiltonian monodromy group of $L = T(x,y)$ is trivial,
\begin{equation}
	\label{eq:mon5}
	\ch_{L}=\{1\}.
\end{equation}
If on the other hand $y$ is irrational and $x \in \Z\langle 1,y \rangle$, then let $k_1,k_2 \in \Z$ be the unique pair with $x = k_1 + k_2 y$. We have
\begin{equation}
	\label{eq:mon4}
	\ch_{L}
	=
	\left\langle
	\begin{pmatrix}
		-1 & 0 \\
		2k_2 & 1
	\end{pmatrix}
	\right\rangle
	=
	\left\{ 
	\begin{pmatrix}
		1 & 0 \\
		0 & 1
	\end{pmatrix},
	\begin{pmatrix}
		-1 & 0 \\
		2k_2 & 1
	\end{pmatrix}
	\right\} \cong \Z_2,
\end{equation}
if least one of the coefficients~$k_1,k_2$ is even. If both $k_1,k_2$ are odd\footnote{This is the case where $(x,y) \in S_y$, see \eqref{eq:badx} in which the billiard construction fails. Therefore we only get an inclusion.}, then
\begin{equation}
	\label{eq:mon4.5}
	\ch_{L}
	\subseteq 
	\left\langle
	\begin{pmatrix}
		-1 & 0 \\
		2k_2 & 1
	\end{pmatrix}
	\right\rangle 
    \cong \Z_2.
\end{equation}

Let $y$ be rational and write $y=\frac{p}{q}$ with $q \in \N_{\geqslant 1}, p \in \Z$ coprime. If $x \notin \frac{1}{q}\Z$, then
\begin{equation}
	\label{eq:mon6}
	\ch_{L} 
	=
	\left\langle
	\begin{pmatrix}
		1 & 0 \\
		2q &  1
	\end{pmatrix}
	\right\rangle
	= 
	\left\{ 
	\left.
	\begin{pmatrix}
		1 & 0 \\
		2kq & 1
	\end{pmatrix}
	\right\vert
	k \in \Z
	\right\} 
	\cong \Z.
\end{equation}

If $x \in \frac{1}{q} \Z$, write $x = \frac{p'}{q}$ and $p' = m_1 q + m_2 p$, where the pair $m_1,m_2 \in \Z$ is chosen such that $m_2 \in \Z_{\geqslant 0}$ is minimal\footnote{Since $p,q$ are coprime, such a pair $m_1,m_2$ exists.}. Then the Hamiltonian monodromy group of $L$ is 

\begin{equation}
	\label{eq:mon7}
	\ch_{L} 
	=
	\left\langle
	\begin{pmatrix}
		1 & 0 \\
		2q &  1
	\end{pmatrix},
	\begin{pmatrix}
		-1 & 0 \\
		2m_2 & 1
	\end{pmatrix}
	\right\rangle
	= 
	\left\{ 
	\left.
	\begin{pmatrix}
		(-1)^{\delta} & 0 \\
		2(\delta m_2 + kq) & 1
	\end{pmatrix}
	\right\vert
	k \in \Z, \delta \in \{0,1\}
	\right\} 
	\cong \Z \rtimes \Z_2,
\end{equation}
if $p' \not\equiv p + q \mod 2$, and it satisfies
\begin{equation}
	\label{eq:mon8}
	\ch_{L} 
	\subseteq
	\left\langle
	\begin{pmatrix}
		1 & 0 \\
		2q &  1
	\end{pmatrix},
	\begin{pmatrix}
		-1 & 0 \\
		2m_2 & 1
	\end{pmatrix}
	\right\rangle,
\end{equation}
if $p' \equiv p + q \mod 2$.
\end{theorem}

\proof
We use the same notation as in the proof of Theorem~\ref{thm:main}. 

First suppose~$(x,y) \in \Sigma_2$. Then~$\cd(x,y) = \{D_1,D_2\}$ and the map~$\phi_* \in \Aut H_2(X_{\alpha},T(x,y))$ either acts by the identity on~$\cd(x,y)$ or permutes its elements by~\cite[Theorem 4.7]{Bre23}. In the former case, the induced monodromy on~$H_1(T(x,y))$ is the identity and in the latter, it is the swap of basis vectors~$e_1 = \pp D_2, e_2 = \pp D_1$. This proves the inclusion~$\subseteq$ in~\eqref{eq:mon2}. To prove equality, consider the symmetric probe with direction vector~$(1,-1)$ passing through~$(x,y)$ together with \eqref{eq:probehomology}.

In the case~$(x,y) = (1,0)$, we have~$\cd(1,0) = \{D_1,D_2,D_3\}$. Since~$\phi$ is a Hamiltonian diffeomorphism it preserves the element~$D_1 + D_3$ which is in the image of the natural map~$H_2(X_{\alpha}) \rightarrow H_2(X_{\alpha}, T(1,0))$. Therefore the only admissible permutations of the elements in~$\cd(1,0)$ are the identity and~$(D_1,D_2,D_3) \mapsto (D_3,D_2,D_1)$. The former induces the identity on~$H_1(T(x,y))$ and the latter the map~$e_1 \mapsto e_1, e_2 \mapsto -e_2$, proving the inclusion~$\subseteq$ in~\eqref{eq:mon3}. To prove equality, consider the vertical probe passing through~$(1,0)$, or a lift to~$S^2 \times S^2$ of an appropriate rotation of the second factor in~$S^2 \times S^2$. The case $(x,y) = (-1,0)$ follows analogously.

Let us now turn to the remaining case where~$(x,y) \in \Sigma_1 \cup \{(0,0)\}$, meaning we prove~\eqref{eq:mon1} and~\eqref{eq:mon1.5}. Note that \cref{lem:mainobstructionmonodromy} does not apply, since $(x,y) \notin Q$. The difference with the case~$(x,y) \in Q$ is that now~$\cd(x,y) = \{D_1,D_3\}$, but the idea of proof is the same. We first treat the case~$\phi_*D_1 = D_1$, for which we obtain \eqref{eq:monodromyQ} by the same arguments. Since~$y = 0$, the equation in~\eqref{eq:monodromyQ2} has no solution for $\delta = 1$, except if~$x = 0$, which is precisely the case \eqref{eq:mon1.5}. Second, suppose that $\phi_*D_1 = D_3$. Similar arguments as in the proof of Lemmata \ref{lem:mainobstruction} and \ref{lem:mainobstructionmonodromy} yield the following cases. For~$\det \phi_* = +1$, we obtain 
\begin{equation}
	\label{eq:mon+1alt}
	(\phi\vert_L)_* 
	=
	\begin{pmatrix}
		-1 & 0 \\
		2l_2 & -1 
	\end{pmatrix}, \quad
	l_1 + l_2y = x, \quad
	l_1, l_2 \in \Z,
\end{equation}
and for~$\det \phi_* = -1$, we obtain
\begin{equation}
	\label{eq:mon-1alt}
	(\phi\vert_L)_*
	=
	\begin{pmatrix}
		1 & 0 \\
		-2k_2 & -1 
	\end{pmatrix}, \quad
	k_1 + k_2y = 0, \quad
	k_1, k_2 \in \Z.
\end{equation}
Again, in $\Sigma_1 \cup \{(0,0)\}$, the second equation in~\eqref{eq:mon+1alt} has no solution except if~$x = 0$. This proves the inclusions~$\subseteq$ in equations~\eqref{eq:mon1} and~\eqref{eq:mon1.5}. All of these monodromy matrices can be realized by noting that every~$(x,y) \in \Sigma_1$ lies on the midpoint of a vertical symmetric probe and a symmetric probe with directional vector~$(1,-1)$. By \eqref{eq:matrices}, these probes yield the elements
$
\begin{pmatrix}
	1 & 0 \\
	0 & -1
\end{pmatrix} \text{ and }
\begin{pmatrix}
	1 & 0 \\
	2 & -1	
\end{pmatrix},
$
respectively. The point~$(0,0)$ additionally lies on the midpoint of a horizontal symmetric probe, yielding the additional generator
$
\begin{pmatrix}
	-1 & 0 \\
	0 & 1
\end{pmatrix}.
$

Let~$(x,y) \in Q$, meaning we can use Lemmata \ref{lem:mainobstructionmonodromy}, \ref{lem:billiardmonodromy} and \cref{cor:monodromyiff} with $x' = x$.

First let~$y \notin \Q$. Then there is no non-trivial pair $k_1,k_2 \in \Z$ which solves~\eqref{eq:monodromyQ2} with $\delta =0$, meaning that there is no element with positive determinant in~$\ch_{L}$. For $\delta = 1$, \eqref{eq:monodromyQ2} becomes $x = k_1 + k_2y$, which has a solution only if $x \in \Z\langle 1,y \rangle$. This proves~\eqref{eq:mon5}. Note that $(x,y) \in S_y \cap Q$ if and only if $x = l_1 + l_2y$ for some $l_1,l_2 \in \Z$ odd. Therefore, \eqref{eq:mon4} follows from \cref{cor:monodromyiff}. This also explains why we do not obtain equality in \eqref{eq:mon4.5}.

Now let~$y = \frac{p}{q} \in \Q$. Again, recall that $(x,y) \in S_y \cap Q$ if and only if $x = l_1 + l_2y$ for some $l_1,l_2 \in \Z$ odd. Note that in that case $x \in \Z\langle 1,\frac{p}{q} \rangle = \frac{1}{q} \Z$. 

Let $x \notin \frac{1}{q} \Z$. Then $(x,y) \notin S_y$. Furthermore, let $\delta,k_1,k_2$ be a triple satisfying \eqref{eq:monodromyQ2} with $x=x'$ in \cref{lem:mainobstructionmonodromy}. Then $\delta = 0$ because in the case $\delta = 1$, \eqref{eq:monodromyQ2} yields $qx = qk_1 + pk_2$ which would imply $x \in \frac{1}{q} \Z$. In the case $\delta = 0$, we find $qk_1 + pk_2 = 0$ and thus~$(k_1,k_2) = k(-p,q)$ for some~$k \in \Z$. This proves \eqref{eq:mon6} by \cref{cor:monodromyiff}.

Let $x = \frac{p'}{q} \in \frac{1}{q}\Z$ for some $p' \in \Z$. Write $p'=m_1q + m_2p$ with $m_2 \in \Z_{\geqslant 0}$ minimal among all $m_2$ for which there is $m_1 \in \Z$ solving this equation. Let $\delta,k_1,k_2$ be a triple satisfying \eqref{eq:monodromyQ2} with $x=x'$. As above, if $\delta = 0$, then $k_2 = kq$ for some $k \in \Z$. If $\delta = 1$, then we obtain the equation $m_1q + m_2p = k_1q + k_2p$ from \eqref{eq:monodromyQ2}. Since $p,q$ are coprime, this implies that $k_2 \equiv m_2 \mod q$, proving \eqref{eq:mon8} and the inclusion $\subseteq$ in \eqref{eq:mon7}. We now need to prove that $(x,y) \in S_y$ if and only if $p' \equiv p + q \mod 2$. Recall from \eqref{eq:badx} that for $(x,y) \in Q$, we have $(x,y) \in S_y$ if and only if there are $l_1,l_2 \in \Z$ such that $x = 2l_1 + 1 + (2l_2 + 1)y$. Setting $y = \frac{p}{q}$ and $x = \frac{p'}{q}$, we find $p' = (2l_1 + 1)q + (2l_2 + 1)p$. Since $p,q$ are coprime there are $l_1,l_2$ satisfying the latter equation if and only if $p' \equiv p + q \mod 2$. This proves \eqref{eq:mon7} and shows that we cannot improve the inclusion \eqref{eq:mon8} by our methods. 
\proofend

\subsection{Cautionary example}
\label{ex:badtorus}

Let $\alpha > 1$ and consider the torus $L = T(0,1) \subset X_{\alpha}$. The methods used in this paper are insufficient to determine the Hamiltonian monodromy group $\ch_{L} \subset \Aut H_1(L)$, as we will discuss now. This illustrates that the hypothesis in \cref{lem:billiardmonodromy} of the billiard trajectory not hitting any corners of the billiard table is necessary. 

In terms of obstructions, note that we are in the case of \eqref{eq:mon8} of \cref{thm:monodromy} with $p=q=1$ and $p'=m_2=0$ and thus 
\begin{equation}
    \label{eq:monodromyT01}
	\ch_{T(0,1)} 
	\subseteq
	\left\langle
	\begin{pmatrix}
		1 & 0 \\
		2 &  1
	\end{pmatrix},
	\begin{pmatrix}
		-1 & 0 \\
		  0 & 1
	\end{pmatrix}
	\right\rangle
	= 
	\left\{ 
	\begin{pmatrix}
		\pm 1 & 0 \\
		2k & 1
	\end{pmatrix}
	\right\} 
	\cong \Z \rtimes \Z_2.
\end{equation}
Alternatively, this special case follows straightforwardly from \cref{lem:mainobstructionmonodromy} by noting that every pair $\delta,k_2$ solves \eqref{eq:monodromyQ2} with $y=1$ and $x=x'=0$ by setting $k_1 = -k_2$. The second generator $e_1 \mapsto -e_1, e_2 \mapsto e_2$ in \eqref{eq:monodromyT01} can be realized by a horizontal probe or lifting to $S^2 \times S^2$ an appropriate rotation of the first factor in $S^2 \times S^2$, in which the circle factor of $T(0,1)$ sits as the equator. We obtain
\begin{equation}
    \label{eq:monodromyT012}
	\Z_2 
    \cong 
    \left\{ 
	\begin{pmatrix}
		\pm 1 & 0 \\
		  0 & 1
	\end{pmatrix}
	\right\} 
    \subseteq
    \ch_{T(0,1)} 
	\subseteq
	\left\{ 
	\begin{pmatrix}
		\pm 1 & 0 \\
		2k & 1
	\end{pmatrix}
	\right\} 
	\cong \Z \rtimes \Z_2,
\end{equation}
which cannot be improved by the methods of this paper. Indeed, the only admissible symmetric probes containing $(0,1)$ are horizontal and vertical. There are no symmetric probes with directional vectors $(\pm 1, -1)$, since both segments hit the vertices of $\square_1$. Equivalently, the good billiard trajectory of $(0,1)$ has one single admissible bouncing point in the sense of \cref{def:goodbilliard}. 

\section{Lagrangian tori from ball embeddings}

\label{sec:ballembeddable}

The main goal of this section is proving the results discussed in \S\ref{ssec:introballemb}. We start by giving some background material on ball embeddings into $X_{\alpha}$ in \S\ref{ssec:conn}. Before proving Theorems \ref{thm:notexotic} and \ref{thm:ballcases} in \S\ref{ssec:proofnotexotic} and \S\ref{ssec:proofballcases}, respectively, we define and discuss the Chekanov torus in \S\ref{ssec:chektorus}.

\subsection{Ball embeddings into $X_{\alpha}$}

\label{ssec:conn}

Recall that we consider symplectic embeddings $\varphi \colon B^4(r)\hookrightarrow X_{\alpha}$ of the standard closed ball of capacity $r$ into $X_{\alpha}$. By Gromov's nonsqueezing theorem \cite{Gro85}, the Gromov width of $X_{\alpha}$ is $2\alpha$ and hence such an embedding exists only if $r < 2\alpha$. See also \cite[Section 5]{Bir97}. A crucial ingredient for the results of this section is the connectedness of the space of ball embeddings into $X_{\alpha}$, which goes back to McDuff \cite{McD98}. The formulation we will use is as follows. 

\begin{lemma}[McDuff]
\label{lem:ballconn}
Let $\varphi,\varphi' \colon B^4(r) \hookrightarrow X_{\alpha}$ be symplectic ball embeddings. Then there is $\psi \in \Ham (X_{\alpha})$ with $\varphi' = \psi \circ \varphi$.
\end{lemma}

This follows from \cite{McD98}. Indeed, by the connectedness of the space of ball embeddings, there is a smooth family of ball embeddings interpolating between $\varphi, \varphi'$. This isotopy can be extended to an isotopy of slightly larger balls, see for example \cite[Theorem 3.3.1]{McDSal17}. Since the ball is simply connected, the symplectic vector field of the isotopy is induced by a globally defined Hamiltonian vector field obtained by performing a cut-off in the larger ball. 

By the connectedness of ball embeddings, we can use a reference ball embedding. Let
\begin{equation}
    \label{eq:phizero}
    \varphi_0 \colon \Int B^4(2\alpha) \hookrightarrow X_{\alpha} 
\end{equation}
be a $T^2$-equivariant symplectic embedding of the open ball at the lower left vertex of the moment map image $\square_{\alpha}$, see also the left-hand side of \cref{fig:4}. Throughout the section, we will use its restrictions to smaller closed balls, $\varphi_0 \colon B^4(r) \hookrightarrow X_{\alpha}$ for all $r < 2\alpha$ as reference symplectic ball embeddings. By abuse of notation, we denote them by $\varphi_0$, as well. By equivariance, the embedding $\varphi_0$ maps product tori $\Theta(b,c) = S^1(b) \times S^1(c)$ in $\R^4$ to split tori in $X_{\alpha}$. More precisely, for $\varphi_0 \colon B^4(r) \hookrightarrow X_{\alpha}$, we find
\begin{equation}
    \label{eq:producttosplit}
    \varphi_0(\Theta(b,c)) = T(b-\alpha-1,c-\alpha)
\end{equation}
for all $b,c > 0$ with $b+c \leqslant r$. 

Recall from \cref{def:types} that a Lagrangian torus $L \subset X^4$ is called \emph{ball-Clifford}, \emph{ball-Chekanov} or \emph{ball-nonmonotone} if it is the image of a \emph{Clifford}, \emph{Chekanov} or \emph{nonmonotone product torus} in $B^4(r) \subset \R^4$ under some ball embedding. 

\begin{lemma}
\label{lem:standardballembedding}
For any Lagrangian torus $L \subset X_{\alpha}$ we have that $L$ is ball-$\tau$ for $\tau$ one of the types (Clifford, Chekanov or nonmonotone) if and only if there is a Lagrangian torus $T \subset B^4(r)$ for some $0<r<2\alpha$ of the same type $\tau$ satisfying $\varphi_0(T) \cong L$ for $\varphi_0$ the reference embedding \eqref{eq:phizero}.
\end{lemma}

\proof 
For the \emph{if}-direction let $\psi \in \Ham (X_{\alpha})$ with $\psi(L) = \varphi_0(T)$. Then the symplectic ball embedding $\psi^{-1} \circ \varphi_0$ proves the claim. For the \emph{only if}-direction, let $L$ be ball-$\tau$. Then by definition, there is a ball embedding $\varphi \colon B^4(r) \hookrightarrow X_{\alpha}$ for which $\varphi^{-1}(L)$ is of type $\tau$. Set $T = \varphi^{-1}(L)$. By \cref{lem:ballconn}, there is $\psi \in \Ham (X_{\alpha})$ such that $\psi \circ \varphi = \varphi_0 \vert_{B^4(r)}$. Recall that $r < 2\alpha$ since the latter is the Gromov width of $X_{\alpha}$. Therefore, we find $\psi(L) = \varphi_0(T)$. 
\proofend

\subsection{The Chekanov torus}
\label{ssec:chektorus}

Let us briefly recall the definition of the Chekanov torus in $\C^2 = \{(z_1,z_2)\}$ via the lifting of curves from a symplectic quotient, which is essentially based on the viewpoint of Eliashberg--Polterovich \cite{EliPol97}. See \cite{Bre22} for a detailed exposition or \cite{Che96,CheSch10,Aur07} for slightly different points of view. The Hamiltonian $H(z_1,z_2) = \pi\vert z_1 \vert^2 - \pi \vert z_2 \vert^2$ generates a Hamiltonian circle action. Its level set $H^{-1}(0)$ is singular, but one can still make sense of symplectic reduction by $H$ in the complement of the singular point $(0,0) \in H^{-1}(0)$. The symplectic quotient carries a residual Hamiltonian $S^1$-action induced by the $T^2$-action on $\C^2$ and it can be equivariantly identified with $\C^{\times}$ equipped with the restriction of the standard symplectic form of $\C$ and the standard Hamiltonian $S^1$-action by rotation. Denote the projection map by $\pi_0 \colon H^{-1}(0)\setminus \{(0,0)\} \rightarrow \C^{\times}$. This map establishes a bijection between Lagrangian tori in $H^{-1}(0)\setminus \{(0,0)\}$ and smooth closed embedded curves in $\C^{\times}$. In particular, the Clifford torus $T_{\rm Cl}(a) = T(a,a)$ corresponds to the standard circle $S^1(a) \subset \C^{\times}$ enclosing symplectic area $a>0$. Since area preserving isotopies of $\C^{\times}$ lift to Hamiltonian isotopies in $\C^2$, every closed curve isotopic to some $S^1(a)$ lifts to a Lagrangian torus which is Hamiltonian isotopic to $T_{\rm Cl}(a)$. On the other hand, it turns out that any curve which is not isotopic to any $S^1(a)$ in $\C^{\times}$ lifts to a copy of the Chekanov torus. 

\begin{definition}
For every $a > 0$, let $\gamma_a \subset \C^{\times}$ be a curve bounding a topological disk of area $a$ in $\C^{\times}$. Then we call its lift $T_{\rm Ch}(a) = \pi_0^{-1}(\gamma_a) \subset \C^2$ a \emph{Chekanov torus}. 
\end{definition}

This yields a monotone Lagrangian torus which is well-defined up to Hamiltonian isotopy. By choosing a suitable curve, one can show that the symplectic ball $B^4(r)$ contains a representative of $T_{\rm Ch}(a)$ for all $a < \frac{r}{2}$. In particular, given a symplectic ball embedding $\varphi \colon B^4(r) \hookrightarrow X$, we can consider the Lagrangian torus $\varphi(T_{\rm Ch}(a)) \subset X$ for every $a < \frac{r}{2}$. 

\subsection{Proof of Theorem \ref{thm:notexotic}} 
\label{ssec:proofnotexotic}

\begin{figure}
  \centering
  \begin{tikzpicture}
    \begin{scope}[scale=1,shift={(-4,0)}]
        \fill[black!5] (-3,-2) rectangle (3,2);
        \fill[black!20] (-3,-2) -- (-3,2) -- (1,-2);
        \draw[black!60, thick, dashed] (-3,2) -- (1,-2);
        \draw[black!60, thick] (-3,-2) rectangle (3,2);
        \draw[very thick,black] (-3,-2)--(1,2);
        \fill[thick, black] (-3,-2)  circle[radius=1.8pt];
        \node at (-0.7,-1.7) {$\varphi_0(\Int B^4(2\alpha))$};
        \node at (2.5,-1.7) {$\square_{\alpha}$};
        \node at (0.65,1.3) {$\sigma$};
        \node at (-3,-2.3) {$\mu(p_{\rm sing})$};
    \end{scope}
    \begin{scope}[scale=0.7,shift={(4,0)}]
        \filldraw[color=black!60, fill=black!5, very thick](0,0) circle (3);
			\fill[thick, black] (0,-3)  circle[radius=2.5pt];
			\draw (0,-3.4) node{$s$};
                \begin{scope}[scale=0.5,shift={(0,-3.15)}]
                    \draw[black, very thick] (1,2) ..controls +(-0.8,0) and + (0.2,0.7).. (-0.5,0.5) ..controls +(-0.2,-0.7) and +(0,1).. (-1.5,-1.5) ..controls +(0,-1) and +(-0.7,-0.7).. (1.5,-1.5) ..controls +(0.7,0.7) and +(1,0).. (1,2);
                \end{scope}
                \begin{scope}[scale=0.55,shift={(0,4.5)}]
                    \draw[black, very thick] (-3,0) .. controls +(0,-0.3) and +(-0.6,0).. (0,-0.4) .. controls +(0.6,0) and +(0,-0.3) .. (3,0);
			        \draw[dashed ,black, very thick] (-3,0) .. controls +(0,0.3) and +(-0.6,0).. (0,0.4) .. controls +(0.6,0) and +(0,0.3) .. (3,0);
                \end{scope}
            \draw[black!50, thick] (-3,0) .. controls +(0,-0.3) and +(-0.6,0).. (0,-0.4) .. controls +(0.6,0) and +(0,-0.3) .. (3,0);
			\draw[dashed,black!50, thick] (-3,0) .. controls +(0,0.3) and +(-0.6,0).. (0,0.4) .. controls +(0.6,0) and +(0,0.3) .. (3,0);
			\node at (1.4,-2) {$\gamma_a$};
            \node at (1.3,1.8) {$C(a)$};
            \node at (3.45,-2) {$S^2(2\alpha)$};
    \end{scope}
  \end{tikzpicture}
  \caption{The main idea of the proof of \cref{thm:notexotic}. Symplectic $S^1$-reduction on $\sigma$ yields the punctured spere on the right-hand side. Any Hamiltonian isotopy mapping $\gamma_a$ to $C(a)$ lifts to a Hamiltonian isotopy in $X_{\alpha}$ which maps $T_{\rm Ch}(a)$ to a split torus.}
  \label{fig:4}
\end{figure}
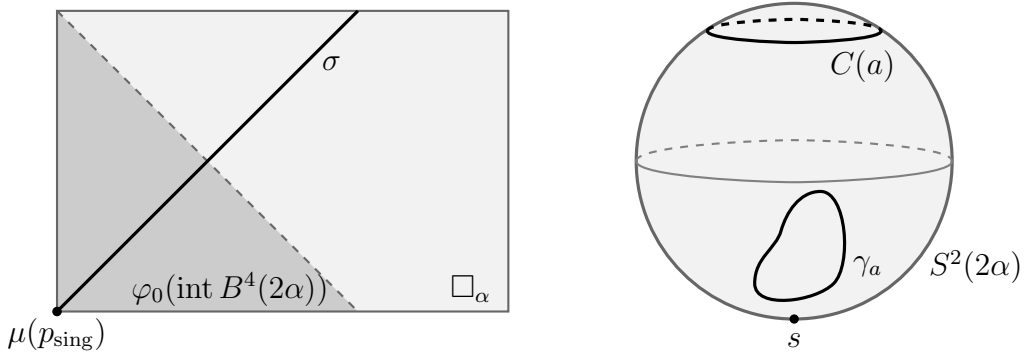

Recall that $r < 2 \alpha$. By \cref{lem:ballconn}, we can assume that $\varphi = \varphi_0$, where $\varphi_0$ is the reference ball embedding \eqref{eq:phizero}. Let $\mu = (\mu_1,\mu_2) \colon X_{\alpha} \rightarrow [-\alpha-1,\alpha+1]$ be the standard moment map on $X_{\alpha}$. The $T^2$-equivariance of $\varphi_0$ translates to $\varphi_0^*(\mu + (\alpha + 1,\alpha)) = (\pi\vert z_1 \vert^2 , \pi \vert z_2 \vert^2)$, where the latter is the standard moment map on $\C^2$. In particular, the Hamiltonian $G = \mu_1 - \mu_2 + 1$ on $X_{\alpha}$ pulls back to $H$ under $\varphi_0$.

We perform symplectic reduction at $G = 0$. The level $\{G=0\}$ contains one singular point, namely $p_{\rm sing} = \mu^{-1}\{(-1-\alpha,-\alpha)\}$ and maps to the segment $\sigma = \{ x - y + 1 = 0 \} \cap \square_{\alpha}$ under $\mu$, see \cref{fig:4}. Crucially, the $S^1$-action generated by $G$ is free elsewhere and in particular on the fibre of $\mu$ over the endpoint $(\alpha-1,\alpha)$ of $\sigma$. This follows from the fact that $\sigma$ intersects the corresponding edge integrally transversely at that point. We deduce that we can perform symplectic reduction on $G^{-1}(0) \setminus \{ p_{\rm sing} \}$. The reduced space carries a residual Hamiltonian $S^1$-action induced by the $T^2$-action on $X_{\alpha}$ and can be equivariantly identified with the punctured two-sphere $S^2(2\alpha) \setminus \{s\}$ (where $s \in S^2$ denotes the south pole) equipped with the standard area form of total area $2\alpha$ and the Hamiltonian $S^1$-action by rotation. Note that under this identification, the fibre $\mu^{-1}(\alpha-1,\alpha)$ maps to the north pole of $S^2$. We denote the projection map by $\pi \colon G^{-1}(0) \setminus \{ p_{\rm sing} \} \rightarrow S^2(2\alpha) \setminus \{s\}$. Since split tori are $T^2$-invariant, they project to the circles of constant height in $S^2(2\alpha) \setminus \{s\}$. For every $0 < b < 2\alpha$, let $C(b)$ denote the circle of constant height bounding a disk of area $b$ in $S^2 \setminus \{s\}$. More precisely, we obtain
\begin{equation}
    \label{eq:projofsplittori}
    \pi(T(\alpha - b - 1, \alpha - b)) = C(b)
\end{equation}
for all $0 < b < 2\alpha$. We now consider the Chekanov torus $T_{\rm Ch}(a) \subset B^4(r)$ for $a < \frac{r}{2} < \alpha$ and its embedding under $\varphi_0$. The reduced space under $H$ of $B^4(r)$ is a standard two-disk of area $\frac{r}{2}$ punctured at the origin. Since $\varphi_0^* G = H$, this punctured disk can be identified with the disk (punctured at the south pole) in the lower hemisphere of $S^2(2\alpha) \setminus \{s\}$ obtained as union of circles of constant height lower than $\frac{r}{2}$. Under this identification, we find that $\pi(T_{\Ch}(a)) = \gamma_a$, where $\gamma_a$ is a curve bounding a disk of area $a$ in $S^2(2\alpha) \setminus \{s\}$. Since they bound disks of the same area, the curve $\gamma_a$ is Hamiltonian isotopic to $C(a)$ in $S^2(2\alpha) \setminus \{s\}$. See the right-hand side of \cref{fig:4} for a sketch of the reduced space and the curves therein. By \eqref{eq:projofsplittori} and the fact that Hamiltonian isotopies can be lifted from symplectic quotients, this proves the claim.
\proofend

\subsection{Proof of Theorem~\ref{thm:ballcases}}
\label{ssec:proofballcases}

Throughout the proof, we use the reference ball embedding $\varphi_0 \colon B^4(r) \hookrightarrow X_{\alpha}$ from \eqref{eq:phizero} and in each case, we choose $r < 2\alpha$ appropriately without stating it explicitly. 

For the ball-Clifford case, note that the image of Clifford tori $T_{\rm Cl}(a) = \Theta(a,a)$ for $a \in (0,\alpha)$ is $\varphi_0(T_{\rm Cl}(a)) = T(a - \alpha - 1, a - \alpha)$ by \eqref{eq:producttosplit}. The corresponding set of points in $\square_{\alpha}$ is $\Sigma \cap \{x < 0,y < 0\}$. By \cref{lem:standardballembedding}, a Lagrangian torus $L$ is ball-Clifford if and only if $L \cong T(a - \alpha - 1, a - \alpha)$ for some $a \in (0,\alpha)$. Applying this to split tori, \cref{thm:main} yields that $T(x,y)$ is ball-Clifford if and only if $(x,y) \in \Sigma$ with $y \neq 0$. This proves (1) in \cref{thm:ballcases}. 

For the ball-Chekanov case, note that the image of Chekanov tori $T_{\rm Ch}(a)$ for $a \in (0,\alpha)$ is $\varphi_0(T_{\rm Ch}(a)) \cong T(\alpha - a - 1, \alpha - a)$ by \cref{thm:notexotic}. By \cref{lem:standardballembedding}, a Lagrangian torus $L$ is ball-Chekanov if and only if $L \cong T(\alpha - a - 1, \alpha - a)$ for some $a \in (0,\alpha)$. For a split torus $T(x,y)$, it thus suffices to show that $T(x,y) \cong T(\alpha - a - 1, \alpha - a)$ for some $a \in (0,\alpha)$ if and only if $(x,y) \notin \Sigma$ and the good billiard trajectory of $(x,y)$ has a corner as bouncing point. For this, note that the base points $(\alpha - a - 1, \alpha - a)$ lie on the diagonal emanating from the lower left corner, meaning that their good billiard trajectory contains the lower left corner as bouncing point. Thus the claim follows from \cref{thm:main} and the fact that we can cover all four corners using the symmetries $(x,y) \mapsto (-x,y)$ and $(x,y) \mapsto (x,-y)$.

For the ball-nonmonotone case, let $T(x,y) \subset X_{\alpha}$ be a split torus. By \eqref{eq:producttosplit} there is a non-monotone (meaning $b \neq c$) product torus $\Theta(b,c)$ with $\varphi_0(\Theta(b,c)) = T(x,y)$ for some $r < 2\alpha$ if and only if $(x,y) \in \cR$ with
\begin{equation}
    \mathcal{R}=\{ (x,y)\in \Int\square_\alpha \mid y<-(x+1),\ y\ne x+1 \}.
\end{equation}
By \cref{lem:standardballembedding}, a Lagrangian torus $L$ is ball-nonmonotone if and only if $L \cong T(x,y)$ for $(x,y) \in \mathcal{R}$. In particular, no split torus with base point in $\Sigma$ is ball-nonmonotone. Let us turn to the complement of $\Sigma$. In order to apply (3) from \cref{thm:main2}, we restrict our attention to $Q$ as defined in \eqref{eq:Qdefset} and consider equivalently $-\cR$ instead of $\cR$. Note that
\begin{equation}
    \label{eq:rcapq}
    (-\cR) \cap Q = \{ (x,y) \in \Int\square_{\alpha} \sth 1 - y < x < 1 + y \}.
\end{equation}
By the above, we find that if $T(x,y)$ is ball-nonmotonotone if and only if it is Hamiltonian isotopic to a split torus with base point in \eqref{eq:rcapq}. By \cref{thm:main2}, this is equivalent to the existence of some $k_1,k_2 \in \Z$ such that 
\begin{equation}
    \label{eq:segmenty}
    1-y < \pm x + 2k_1 + 2k_2y < 1+y.
\end{equation}
For every $y \in \R \setminus \Q$, such $k_1,k_2$ do exist, since the subgroup $\Z\langle 2, 2y \rangle \subset \R$ is dense in that case. Suppose now that $y = \frac{p}{q} \in \Q$ for $p,q$ coprime. Note that in that case $\Z \langle 2, 2y\rangle = \frac{2}{q} \Z$. Since the interval $(1-y,1+y)$ has length $2y = \frac{2p}{q}$, the only possibility for $(1-y,1+y)$ to not contain values of the form \eqref{eq:segmenty} is that $p = 1$ and $\pm x$ lies on its boundary modulo $\frac{2}{q}$. This is the case if and only if $(x,y) \in \cd \cap \{y > 0\}$. The claim for split tori which are not contained in $Q$ follows from applying a symmetric probe mapping that point into $Q$. \proofend

\proofof{\cref{cor:cheknonmonotone}} By \cref{thm:notexotic}, we have $\varphi(T_{\rm Ch}(a)) \cong T(\alpha - a -1, \alpha - a)$. By \cref{thm:ballcases} (1) and (3), there is a ball embedding $\varphi' \colon B^4(r') \hookrightarrow X_{\alpha}$ and a product torus $\Theta(b,c) \subset B^4(r')$ with $\varphi'(\Theta(b,c)) = T(\alpha - a - 1, \alpha - a)$ if and only if $a \notin \left\{ \alpha - \frac{1}{k} \sth k \in \N_{\geqslant 1} \right\}$. By \cref{lem:ballconn} we find $\varphi'(\Theta(b,c)) \cong \varphi(\Theta(b,c))$ and thus the claim follows. \proofend

\section{Lagrangian packing}

\label{sec:Lagpacking}

In this section, we prove and discuss the results from \S\ref{ssec:introLagpacking}. For the reader's convenience, we start by proving that \cref{thm:MSPS} is a reformulation of \cite[Theorem C]{PolShe23} by adapting its statement to our notational conventions. 

\subsection{Discussion of \cref{thm:MSPS}}

\label{ssec:discussionMSPS}

In \cite{PolShe23}, $S^2 \times S^2$ is equipped with the symplectic form $\omega_{S^2} \oplus 2a \omega_{S^2}$ for $0 < 2a < 1$ and the torus they consider is $L = S^1(B) \times S^1_{\rm eq}$, where $S^1(B)$ is an embedded circle bounding a disk of area $B > 0$ and $S^1_{\rm eq}$ is the equator. For any $k \in \N_{\geqslant 2}$, it follows from \cite[\S 3.4-3.6]{PolShe23} that $\#_P(L) \leqslant k$ whenever 
\begin{equation}
\label{eq:aBC}
    B - C > a, \quad 
    \text{ where } \quad
    C = \frac{1 - 2B}{k-1}.
\end{equation}
Switching back to our notation, let $\alpha > 0$ and $T(x,0) \subset X_{\alpha}$ for $x \in (-1,1)$. By scaling the symplectic form $\omega_{\alpha}$ of $X_{\alpha}$ by a factor of $\frac{1}{2(1 + \alpha)}$ and substituting 
\begin{equation}
    a = \frac{\alpha}{2(1 + \alpha)}, \quad
    B = \frac{1}{2} - \frac{\vert x \vert}{2(1 + \alpha)},
\end{equation}
we see that \eqref{eq:aBC} implies \eqref{eq:PSmain}. 

\subsection{Packing by split tori}

\begin{definition}
For $(x,y) \in \R^2 \setminus \Sigma$, let $b(x,y) \in \N \cup \{\infty\}$ be the cardinality of the set of admissible bouncing points of the good billiard trajectory of at least one of the points $(\pm x, \pm y)$.
\end{definition}

\cref{thm:main} immediately implies the following. 

\begin{corollary}
\label{cor:lowerbound}
Let $L = T(x,y) \subset X_{\alpha}$ be a split torus. Then,
    \begin{equation}
        b(x,y) = \#_{TP}(L) \leqslant \#_{P}(L). 
    \end{equation}
\end{corollary}

Here, $\#_{TP}(L)$ denotes the packing number \emph{by disjoint split tori} of $X_{\alpha}$. Note that it makes sense for general toric symplectic manifolds to ask what the packing number \emph{by disjoint toric fibres} is, see the discussion surrounding \cref{q:packing2}. Using the same notation as in \S\ref{ssec:monodromythm}, we can determine $b(x,y)$ for all $(x,y) \in \R^2$. The distinction of cases follows a similar pattern as in \cref{thm:monodromy}.

\begin{proposition}
\label{prop:toricpacking}
Let $L = T(x,y)$ by a split torus. If $(x,y) \in \Sigma$ then 
\begin{equation}
    \label{eq:bfirst}
    b(x,y) =
    \begin{cases}
        1 & \text{if } (x,y)=(0,0), \\
        2 & \text{if } (x,y) \in \Sigma_1 \cup \{(\pm 1, 0)\}, \\
        4 & \text{if } (x,y) \in \Sigma_2.
    \end{cases}
\end{equation}
If $(x,y) \notin \Sigma$, we suppose without loss of generality that $(x,y) \in Q$. In that case, if $y \in \R \setminus \Q$, then
\begin{equation}
    \label{eq:bsecond}
    b(x,y) = \infty.
\end{equation}
If $y \in \Q$, we write $y = \frac{p}{q}$ for $p,q$ coprime and have
\begin{equation}
    \label{eq:bthird}
	b(x,y) = 
	\begin{cases}
		4p + 2q - 4, & x = \frac{p'}{q} \text{ for some } p'\in \Z \text{ with } p' \equiv p+q \mod 2,\\
		4p + 2q, & x = \frac{p'}{q} \text{ for some } p'\in \Z \text{ with } p' \not\equiv p+q \mod 2,\\
		8p + 4q, & \text{else.}
	\end{cases}
\end{equation}
\end{proposition}

\proof 
The cases in \eqref{eq:bfirst} are an immediate consequence of our definition of admissible bouncing points of a good billiard trajectory in the case $(x,y) \in \Sigma$, see \cref{def:goodbilliard}. The case \eqref{eq:bsecond} follows from the elementary fact that good billiard trajectories on rectangular billiard tables are not periodic and thus have infinitely many bouncing points. Note that such a good billiard trajectory hits a corner of the billiard table \emph{at most once} and hence always has infinitely many admissible bouncing points. Let us turn to the case \eqref{eq:bthird}. Recall from \cref{thm:main2} that points $(x,y),(x',y) \in Q$ correspond to equivalent tori if and only if $x' \pm x \in \Z \langle 2,2y \rangle$. Here, we have $\Z \langle 2,2y \rangle = \frac{2}{q}\Z$. This means that the orbit in $Q \cap \square_y$ is generated by $x \mapsto -x$ and translations by $\frac{2}{q}$. In other words, the segment $[ 0, \frac{1}{q}]$ contains exactly one admissible bouncing point of every good billiard trajectory in $\square_y$. The endpoints $0$ and $\frac{1}{q}$ correspond to the first two cases in \eqref{eq:bthird} depending on the parity of $p+q$ and the points in the interior of $[0,\frac{1}{q}]$ to the last case in \eqref{eq:bthird}. The perimeter of $\square_{y}$ for $y = \frac{p}{q}$ is $4\left( \frac{2p +q}{q} \right)$, which implies the last equality in \eqref{eq:bthird} by dividing by the length $\frac{1}{q}$ of $[0,\frac{1}{q}]$. For the other two equalities we a priori take half of $4(2p + q)$ because the corresponding representatives in $[ 0, \frac{1}{q}]$ are boundary points. Additionally, we need to take the case into account where the good billiard trajectory hits the corner points of the billiard table. By \eqref{eq:badx}, as in the proof of \cref{thm:monodromy}, this happens if and only if $p' \equiv p + q$ in which case we remove four bouncing points to obtain the count of admissible bouncing points. 
\proofend

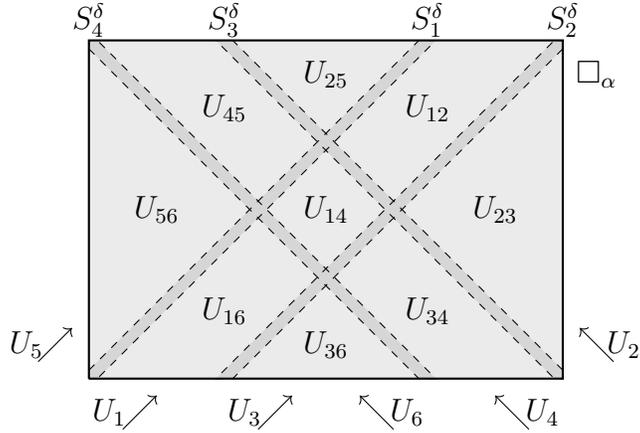
\begin{figure}
		\begin{tikzpicture}[scale=0.45]	
			\fill[black!8] (-7,-5)--(7,-5)--(7,5)--(-7,5)--(-7,-5);
			
			\fill [black!16] (-7,-4.75)--(2.75,5)--(3.25,5)--(-6.75,-5)--(-7,-5)--(-7,-4.75);
			\fill [black!16] (7,4.75)--(-2.75,-5)--(-3.25,-5)--(6.75,5)--(7,5)--(7,4.75); 
			\fill [black!16] (7,-4.75)--(-2.75,5)--(-3.25,5)--(6.75,-5)--(7,-5)--(7,-4.75);
			\fill [black!16] (-7,4.75)--(2.75,-5)--(3.25,-5)--(-6.75,5)--(-7,5)--(-7,4.75); 
			
			\draw [dashed] (-7,-4.75)--(2.75,5)--(3.25,5)--(-6.75,-5)--(-7,-5)--(-7,-4.75);
			\draw [dashed] (7,4.75)--(-2.75,-5)--(-3.25,-5)--(6.75,5)--(7,5)--(7,4.75); 
			\draw [dashed] (7,-4.75)--(-2.75,5)--(-3.25,5)--(6.75,-5)--(7,-5)--(7,-4.75);
			\draw [dashed] (-7,4.75)--(2.75,-5)--(3.25,-5)--(-6.75,5)--(-7,5)--(-7,4.75);
			
			\draw [thick] (-7,-5)--(7,-5)--(7,5)--(-7,5)--(-7,-5); 
			\node at (0,0){$U_{14}$};
			\node at (3,3){$U_{12}$};
			\node at (5,0){$U_{23}$};
			\node at (3,-3){$U_{34}$};
			\node at (0,-4){$U_{36}$};
			\node at (-3,-3){$U_{16}$};
			\node at (-5,0){$U_{56}$};
			\node at (-3,3){$U_{45}$};
			\node at (0,4){$U_{25}$};
			
			\node at (8,4){$\square_{\alpha}$};
			
			\draw [black,->] (2,-6.5)--(1,-5.5);
			\draw [black,->] (6,-6.5)--(5,-5.5);
			\draw [black,->] (8.5,-4.5)--(7.5,-3.5);
			\node at (2.4,-6){$U_6$};
			\node at (6.4,-6){$U_4$};
			\node at (8.8,-4){$U_2$};
			
			\draw [black,->] (-2,-6.5)--(-1,-5.5);
			\draw [black,->] (-6,-6.5)--(-5,-5.5);
			\draw [black,->] (-8.5,-4.5)--(-7.5,-3.5);
			\node at (-2.4,-6){$U_3$};
			\node at (-6.4,-6){$U_1$};
			\node at (-8.85,-4){$U_5$};
			
			\node at (3,5.6){$S_1^{\delta}$};
			\node at (7,5.6){$S_2^{\delta}$};
			\node at (-3,5.6){$S_3^{\delta}$};
			\node at (-7,5.6){$S_4^{\delta}$};
		\end{tikzpicture}
	\caption{The subsets $U_i \subset \square_{\alpha}$ with dashed boundaries and their intersections $U_{ij}$. Arrows indicate the direction of the billiard trajectories realized by $\psi_{\delta}$ in the respective $U_i$. The cut-off region $P_{\delta}$ is dark grey.}
	\label{fig:Ui}
\end{figure}

\section{Late Lagrangian recurrence}
\label{sec:latereturn}

In this section, we provide proofs of \cref{prop:latereturnweak} and \cref{thm:latereturn}.\smallskip

\proofof{\cref{prop:latereturnweak}}
Let $L=L_1,L_2,L_3,\ldots$ be the split tori which have as footpoints the admissible bouncing points of the good billiard trajectory of $(x,y)$ with $L = T(x,y)$. Since $r(x,y) \in \R \setminus \Q$, we can assume that $L_i \neq L_j$ for $i \neq j$. Indeed, if the billiard trajectory hits a corner point of $\square_{r}$ in one direction, then it does not hit one in the other direction. Now let $N \in \N$ and define $\psi \in \Ham (X_{\alpha})$ as the Hamiltonian diffeomorphism obtained by concatenation of the first $N-1$ symmetric probes corresponding to that billiard trajectory with supports chosen small enough so that the supports of the Hamiltonians of probes intersect if and only if the probes intersect. Since the Hamiltonian diffeomorphism associated with the $i$-th probe swaps $L_{i}$ and $L_{i+1}$, we obtain 
\begin{equation}
    \psi(L_1) = L_N, \quad \psi(L_i) = L_{i-1} \quad \forall i \in \{2,\ldots,N\}.
\end{equation}
Since the $L_i$ are pairwise disjoint, this proves the claim.
\proofend

We now turn to the proof of \cref{thm:latereturn}, which requires some preparations. Let $\alpha > 0, \delta > 0$ and define the following line segments in $\square_{\alpha}$
\begin{eqnarray*}
	S_1 = \{ y - x - 1 = 0 \} \cap \square_{\alpha}, \quad
	S_2 = \{ y - x + 1 = 0 \} \cap \square_{\alpha}, \\
	S_3 = \{ y + x - 1 = 0 \} \cap \square_{\alpha} , \quad
	S_4 = \{ y + x + 1 = 0 \} \cap \square_{\alpha}.
\end{eqnarray*}

For every $i$, denote by $S_i^{\delta}$ the set of points in $\square_{\alpha}$ which are at integral affine distance~$< \delta$ from $S_i$. The complement of $S_1^{\delta} \cup S_2^{\delta}$ consists of three connected components which we denote by $U_5,U_1,U_3$ from left to right. Similarly, we denote the three connected components of the complement of $S_3^{\delta} \cup S_4^{\delta}$ by $U_6, U_4, U_2$ from left to right. See also Figure~\ref{fig:Ui}. For every $i$, let $U_i^{\delta}$ be the set of points in $\square_{\alpha}$ which are at integral affine distance $< \delta$ from $U_i$.

For every $i \in \{1,\ldots,6\}$, we define a Hamiltonian diffeomorphism $\varphi_i^{\delta}$ which has support in $\mu^{-1}(U_i^{\delta})$ and which realizes a part of a billiard trajectory for all toric fibres in $\mu^{-1}(U_i)$. Heuristically, every $\varphi_i^{\delta}$ comes from a certain symmetric probe which has been thickened to have as large a support as is permitted by $\square_{\alpha}$ with a cut-off region of width $\delta$. For the reader's convenience, we give an explicit definition of the $\varphi_i^{\delta}$. 

We introduce model spaces $M_1, M_2$ and Hamiltonian diffeomorphisms on them. 
Let $M_1= Z \times S^2$, where $Z$ denotes the cylinder $Z = \{(\theta,p) \in T^*S^1 \, \vert \, -1 < p < 1 \}$ equipped with the product symplectic form $\omega_1 = d\theta \wedge dp \oplus 2\alpha \omega_{S^2}$. Denoting the height function generating the standard Hamiltonian circle action on $(S^2,2\alpha\omega_{S^2})$ by $h \colon S^2 \rightarrow [-\alpha,\alpha]$, the moment map $\nu_1(\theta,p,x)=(p,h(x))$ defines a toric structure on $M_1$. Let $R \in \Ham(S^2,2\alpha\omega_{S^2})$ with $h \circ R = -h$; this is for example a suitable solid rotation. By letting $R$ act on the second factor of $M = Z \times S^2$ and choosing an appropriate cut-off along the $p$-direction, we find $\chi_1^{\delta} \in \Ham(M_1,\omega_1)$ satisfying 
\begin{equation}
	\label{eq:chi_1delta}
	\chi^{\delta}_1(\theta,p,x) = (\theta,p,R(x)), \quad
	\text{for } -1 + \delta < p < 1 - \delta
\end{equation} 
and equal to the identity outside of a slightly larger compact set. The moment map image of $\nu_1$ is $(-1,1) \times [-\alpha,\alpha]$. Note that by \eqref{eq:chi_1delta}, $\chi_1^{\delta}$ preserves the set of toric fibres over $V_1 = (-1+\delta,1-\delta) \times [-\alpha,\alpha] \subset \im(\nu_1)$, where it acts by the map $(x,y) \mapsto (x,-y)$.

The second model space $M_2$ is the standard open symplectic ball $B^4(2\alpha) = \{\pi\vert z_1 \vert^2 + \pi \vert z_2 \vert^2 \leqslant 2\alpha \} \subset \C^2$ equipped with the toric moment map $\nu_2(z_1,z_2)= (\pi \vert z_1 \vert^2, \pi \vert z_2 \vert^2)$. By a cut-off, we can find a Hamiltonian diffeomorphism $\chi_2^{\delta} \in \Ham(B^4(2\alpha))$ satisfying 
\begin{equation}
	\label{eq:chi_2delta}
	\chi^{\delta}_2(z_1,z_2) = (z_2,z_1), \quad \text{for } \pi\vert z_1 \vert^2 + \pi \vert z_2 \vert^2 \leqslant 2\alpha - \delta
\end{equation} 
and equal to the identity outside a slightly larger compact set. By~\eqref{eq:chi_2delta}, $\chi_2^{\delta}$ preserves the set of toric fibres over $V_2 = \{ (x,y) \in \R^2_{\geqslant 0} \, \vert \, x + y \leqslant 2\alpha - \delta \}$, where it acts by the map $(x,y) \mapsto (y,x)$.

For $i \in \{1,4\}$, the pair $(U_i^{\delta},U_i)$ is integral affine equivalent to $(\im \nu_1 , V_1)$, by which we mean that there is a $T^2$-equivariant symplectic embedding $f_i \colon M_1 \hookrightarrow X_\alpha$ with image $\mu^{-1}(U_i^{\delta})$ and mapping $\nu_1^{-1}(V_1)$ to $\mu^{-1}(U_i)$. Set $\varphi_i^{\delta} = f_i \circ \chi_1^{\delta} \circ f_i^{-1}$ on $\mu^{-1}(U_i^{\delta})$ and extend it by the identity to obtain $\varphi_i^{\delta} \in \Ham(X_\alpha)$ for $i \in \{1,4\}$.

For $i \in \{2,3,5,6\}$, we proceed similarly. The pair $(U_i^{\delta},U_i)$ is integral affine equivalent to $(\im \nu_2 , V_2)$, meaning that there is a $T^2$-equivariant symplectic embedding $f_i \colon M_2 \hookrightarrow X_\alpha$ with image $\mu^{-1}(U_i^{\delta})$ and mapping $\nu_2^{-1}(V_2)$ to $\mu^{-1}(U_i)$. Set $\varphi_i^{\delta} = f_i \circ \chi_2^{\delta} \circ f_i^{-1}$ on $\mu^{-1}(U_i^{\delta})$ and extend it by the identity to obtain $\varphi_i^{\delta} \in \Ham(X_\alpha)$ for $i \in \{2,3,5,6\}$.

\begin{proposition}
Let $L = T(x) \in \mathcal{T}$ be a toric fibre with $x \in U_i$, then $\varphi^{\delta}_i$ induces a billiard move on $x$ in the direction $v_i$ for $v_1=v_3=v_5=(1,1)$ and $v_2 = v_4 = v_6 = (1,-1)$. 
\end{proposition}

\proof 
This is an immediate consequence of~\eqref{eq:chi_1delta},~\eqref{eq:chi_2delta} and the fact that the maps $f_i$ are equivariant and thus map toric fibres to toric fibres. 
\proofend

For all $1 \leqslant i < j \leqslant 6$, let $U_{ij} = U_i \cap U_j$. There are nine pairs $i<j$ (or seven for $\alpha < 1$) for which $U_{ij} \neq \varnothing$, see also Figure~\ref{fig:Ui}. Note that $\square_\alpha \setminus P_{\delta} = \cup_{1\leqslant i<j \leqslant 6} U_{ij}$, where $P_{\delta} = S_1^{\delta} \cup S_2^{\delta} \cup S_3^{\delta} \cup S_4^{\delta}$. The following properties of the maps $\varphi_{i}^{\delta}$ follow from the construction, again we refer to Figure \ref{fig:Ui}.

\begin{proposition}
\label{prop:psiaux}
\begin{itemize}
\item[\rm (1)] $\varphi^{\delta}_k\vert_{U_{ij}} = \id$ if $k \notin \{i,j\}$,
\item[\rm (2)] $\varphi^{\delta}_k\vert_{U_{ij}}$ acts by the billiard move induced by $\varphi_i^{\delta}$ if $k=i$ and by the billiard move induced by $\varphi_j^{\delta}$ if $k=j$,
\item[\rm (3)] the compositions $\varphi^{\delta}_i \circ \varphi^{\delta}_j$ and $\varphi^{\delta}_j \circ \varphi_i^{\delta}$ induce a billiard move which has a bouncing point in $U_{ij}$.
\end{itemize}
\end{proposition}

Indeed, let $x \in U_i$, then $\varphi_i$ maps $T(x)$ to some $T(x')$, where $x'$ is one of the two adjacent bouncing points of the billiard trajectory of $x$. This follows immediately from the discussions surrounding \eqref{eq:chi_1delta}, \eqref{eq:chi_2delta} and the definition of $\varphi_i^{\delta}$. 

\begin{definition}
The {\rm $\delta$-billiard map} of $X_{\alpha}$ is defined as the composition
\begin{equation} 
    \label{eq:psidelta}
	\psi_{\delta} = \varphi_6^{\delta} \circ
	\varphi_5^{\delta} \circ
	\varphi_4^{\delta} \circ
	\varphi_3^{\delta} \circ
	\varphi_2^{\delta} \circ
	\varphi_1^{\delta} \in \Ham(X_{\alpha}).
\end{equation} 
\end{definition}

The name is warranted by the fact that outside of the union $P_{\delta} = S_1^{\delta} \cup S_2^{\delta} \cup S_3^{\delta} \cup S_4^{\delta}$ of cut-off regions of the $\varphi_i^{\delta}$, the $\delta$-billiard map takes toric fibres to toric fibres on a billiard trajectory.

\begin{lemma}
\label{lem:psibilliard}
Let $L \in \mathcal{T} = \{\text{toric fibres in } X_{\alpha}\}$ be such that $\psi_{\delta}^k(L) \subset X_{\alpha} \setminus \mu^{-1}(P_{\delta})$ for all $k \in \{1,\ldots,N\}$. Then there exists a finite sequence $\{i_k\}_{1 \leqslant k \leqslant N+1}$ which is strictly increasing or strictly decreasing with
\begin{equation}
	\label{eq:psibilliard}
	\psi_{\delta}^k(L) = L_{(i_k)}, \quad
	\text{for all } k \in \{1,\ldots,N+1\}.
\end{equation}
\end{lemma}

\proof 
We suppose that $L = T(x)$ for $x \notin \Sigma$, otherwise both sides of~\eqref{eq:psibilliard} are independent of $k$ and the claim holds trivially. Recall from \eqref{eq:psidelta} that $\psi_{\delta}$ is a composition of maps $\varphi_i^{\delta}$ for $i \in \{1,\ldots,6\}$. Let $\{M_k\}_{0 \leqslant k \leqslant 6(N+1)}$ be defined by
\begin{equation}
	M_0 = L, \quad 
	M_{k+1} = \varphi^{\delta}_{m_k}(M_k), 
\end{equation}
where $m_k -1$ is the remainder of the division of $k$ by $6$. In other words, $\{M_k\}$ is the sequence of Lagrangian tori obtained by consecutively and cyclically applying the maps $\varphi_i^{\delta}$. Note that $M_{6l} = \psi_{\delta}^{l}(L)$ for all $l \leqslant N$. Since, by hypothesis, the sequence $\psi_{\delta}^{l}(L)$ is contained in $X_{\alpha} \setminus \mu^{-1}(P_{\delta})$ and thus consists of toric fibres by Proposition~\ref{prop:psiaux}, it follows that the same holds for $\{M_k\}$. Let $\{x_k\}$ be the corresponding sequence of base points, i.e.\ $M_k = T(x_k)$ for all $0 \leqslant k \leqslant 6(N+1)$.

The sequence $\{M_k\}$ may contain consecutive Lagrangians which coincide. Let $k_{\rm in}$ be the smallest index for which $M_{k_{\rm in}} \neq L$ and let $k_{\rm fin}$ be the largest index for which $M_{k_{\rm fin}} \neq M_{6N + 6}$. Note that $k_{\rm in} \leqslant 6$ and $k_{\rm fin} \geqslant 6N + 1$ since $\Delta_\alpha \setminus P_{\delta} = \cup_{1\leqslant i<j \leqslant 6} U_{ij}$. For every $k_{\rm in} \leqslant k \leqslant k_{\rm fin}$, the following are well-defined
\begin{equation}
	\label{eq:kstars}
	k_* = \max\{k' < k \sth M_{k'} \neq M_{k}\}, \quad
	k^* = \min\{k' > k \sth M_{k'} \neq M_{k}\}.
\end{equation} 
We show that for all such $k$, the triple $x_{k_*},x_k,x_{k^*}$ of corresponding base points realizes a billiard bounce in $\square_{\alpha}$. Since the subsequences obtained for $k < k_{\rm in}$ and for $k > k_{\rm fin}$ are constant, this proves that the full sequence $\{M_k\}$ traces out a (forward or backward) billiard trajectory of $L$; with possible repetitions of the form $M_k = M_{k+1}$. Since $\{\psi^k_{\delta}(L)\}_{0 \leqslant k \leqslant N}$ is a subsequence of $\{M_k\}$, this proves the claim.

To prove that the triple $x_{k_*},x_k,x_{k^*}$ is obtained by performing a billiard bounce at $x_k$, note that $x_k$ is contained in a unique $U_{ij}$. It follows that
\begin{equation}
	\label{eq:tripleMk}
	M_k = \varphi^{\delta}_{i_1}(M_{k_*}), \quad
	M_{k^*} = \varphi^{\delta}_{i_2}(M_k).
\end{equation}
for some $i_1, i_2 \in \{i,j\}$ with $i_1 \neq i_2$. Indeed, by Proposition~\ref{prop:psiaux} every $\varphi_k^{\delta}$ with $k \notin \{i,j\}$ equals the identity on $U_{ij}$ and therefore $i_1,i_2 \in \{i,j\}$. The claim $i_1 \neq i_2$ follows from the fact that the maps $\varphi_i^{\delta}$ are applied cyclically and by increasing order of $i$. This implies that the sequence $i_k$ in~\eqref{eq:psibilliard} is strictly increasing or strictly decreasing, where the former corresponds to the case $i = i_1 < i_2 = j$ and the latter to the case $i = i_2 < i_1 = j$. \proofend

\proofof{\cref{thm:latereturn}} By Lemma~\ref{lem:psibilliard}, under iterates of $\psi_{\delta}$, all product tori follow their billiard trajectory for as long as this trajectory is disjoint from $P_{\delta}$ and this yields disjoint copies of the tori in question. Thus find that 
\begin{equation}
	\mathcal{T}_{\psi_{\delta},N} \supset \Int \square_{\alpha} \setminus \bigcup_{k=0}^N(P_{\delta})_{(-k)},
\end{equation}
where $\mathcal{T}_{\psi_{\delta},N}$ is defined as in \eqref{eq:mathcaltdef}, and with $(P_{\delta})_{(-k)} = \{(x,y) \sth (x,y)_{(k)} \in P_{\delta}\}$ where $(x,y)_{(k)}$ denotes the $k$-th admissible bouncing point on the trajectory of $(x,y)$, as before. Recall that by $\lambda$, we denote the Lebesgue measure on $\square_{\alpha}$ with normalization $\lambda(\square_{\alpha}) = 1$. Furthermore, for every point in $P_{\delta}$, there is exactly one point in $(P_{\delta})_{(-k)}$ from which it follows that (as the intersection may not be disjoint), 
\begin{equation}
	\lambda\left( \bigcup_{k=0}^N(P_{\delta})_{(-k)} \right) \leqslant (N+1) \vol(P_{\delta}). 
\end{equation}
Since $\vol(P_{\delta}) \rightarrow 0$ for $\delta \rightarrow 0$, the claim follows. \proofend

\section{Open questions and discussion}
\label{sec:questions}

In this section we discuss open questions surrounding our results. 

\subsection{Classification of toric fibres and their Hamiltonian monodromy} After $S^2 \times S^2$, the next natural case in which one can ask for a classification of toric fibres is the one-fold blow up of $\CP^2$. Up to scaling, it supports a one-parametric family of symplectic forms depending on the size of the blow-up. Even in the monotone case, the classification is unknown. 

\begin{question}
\label{q:class1}
What is the classification of toric fibres in the one-fold blow-up of $\CP^2$?
\end{question}

As discussed in \cref{rk:Hirzebruch}, a classification for all symplectic forms would, together with our results, yield a classification of toric fibres in all Hirzebruch surfaces. Another question raised by our results is determining the Hamiltonian monodromy group in the cases which remain open in \cref{thm:monodromy}. 

\begin{question}
\label{q:mon1}
What is the Hamiltonian monodromy group in the cases \eqref{eq:mon4.5} and \eqref{eq:mon8}?
\end{question}

\subsection{Lagrangian packing}
\label{ssec:qpacking}
 
\begin{question}
\label{q:packing1}
    What is the packing number of a given split torus in $X_{\alpha}$? 
\end{question}

The only cases where it is known are as follows 
\begin{equation}
    \label{eq:packingnumbers}
    \#_P(T(x,y))
    =
    \begin{cases}
        1 & \text{ if } (x,y) = (0,0),\\
        2 & \text{ if } y = 0 \text{ and } \vert x \vert \in \left( 0, \frac{1}{3} \right), \\
        k & \text{ if } y = 0 \text{ and } \vert x \vert \in \left( \frac{k-2}{k} ,  \frac{k-1}{k+1} \right) \text{ and }\alpha \text{ small enough,} \\
        \infty & \text{ if } r(x,y) \in \R \setminus \Q.
    \end{cases}
\end{equation}

The first case is simply non-displaceability which follows for example from \cite{EntPol06}, since all other split tori are displaceable. The following two cases are consequences of \cref{thm:MSPS} and maximal packings obtained by packing the larger sphere factor by circles (this is where the smallness of $\alpha$ is required), and the last case is \cref{cor:infpacking}. The case of the split tori $T(\pm 1,0) \subset X_{\alpha}$ is particularly interesting, since it is the boundary case between the rigid region of \cref{thm:MSPS} and the flexible region where our methods apply. Even when the packing number is infinite, one can ask whether a packing configuration $\mathcal{P} \subset X_{\alpha}$ of disjoint copies of a Lagrangian $L$ yielding the infinite packing number is maximal in the sense that
\begin{equation}
    \phi(L) \cap \mathcal{P} \neq \varnothing \quad \text{for all} \quad \phi \in \Ham(X_{\alpha}).
\end{equation}
In many cases, packing configurations by split tori obtained by our methods are not maximal, since more copies which are not split tori can be packed into the larger one of the sphere factors. 

\begin{question}
\label{q:packing2}
    Let $X$ be a toric symplecic manifold. For which toric fibres is the toric packing number $\#_{TP}(L)$ equal to the actual Lagrangian packing number $\#_{P}(L)$?
\end{question}

Note that a classification of toric fibres up to Hamiltonian diffeomorphism always yields the toric packing number, see \cite{Bre23} for more examples of such classifications.

\subsection{Ball-embeddable Lagrangian tori} \label{ssec:qballembeddable} Let $(X,\omega)$ be a four-dimensional symplectic manifold. Denote the set of ball-Clifford, ball-Chekanov and ball-nonmonotone Lagrangian tori in $X$ by $\cl_{\rm Cl}, \cl_{\rm Ch}$ and $\cl_{\rm nm}$, respectively. Recall that for $X=X_{\alpha}$ we have
\begin{equation}
    \cl_{\rm Cl} \cap \cl_{\rm Ch} = \varnothing, \quad
    \cl_{\rm Cl} \cap \cl_{\rm nm} = \varnothing, \quad
    \cl_{\rm Ch} \cap \cl_{\rm nm} \neq \varnothing,
\end{equation}
by \cref{thm:ballcases}. By similar arguments, one finds that these three sets are mutually disjoint for monotone $S^2 \times S^2$. For $X=\CP^2$ the situation is
\begin{equation}
    \cl_{\rm Cl} \cap \cl_{\rm Ch} = \varnothing, \quad
    \cl_{\rm Cl} \cap \cl_{\rm nm} \neq \varnothing, \quad
    \cl_{\rm Ch} \cap \cl_{\rm nm} = \varnothing.
\end{equation}

This suggests the following question.

\begin{question}
\label{q:ballembeddable1}
Are the sets $\cl_{\rm Cl}$ and $\cl_{\rm Ch}$ disjoint for all $X$?
\end{question}

It is crucial that the embedded domains are balls. For polydisks for example, some images of Chekanov tori are Hamiltonian isotopic to images of Clifford tori. This happens for example in $X_{\alpha}$ by arguments similar to those used in the proof of \cref{thm:notexotic}.

We can positively answer \cref{q:ballembeddable1} for tori which are small relative to the size of the balls in which they are contained. 

\begin{proposition}
Let $L=\varphi(T_{\rm Cl}(a))$ for a symplectic ball embedding $\varphi \colon B^4(r) \hookrightarrow X$ and $L'=\varphi'(T_{\rm Ch}(a'))$ for a symplectic ball embedding $\varphi' \colon B^4(r') \hookrightarrow X$ with 
\begin{equation}
    \label{eq:relsizes}
    a,a' < \frac{1}{3}\min\{r,r'\}, 
\end{equation}
then $L \not\cong L'$.
\end{proposition}

This follows from a local computation of displacement energy germs as in \cite{Bre23b,BreHauSch23}. The inequalities \eqref{eq:relsizes} are used in certain area-estimates, due to \cite{CheSch16}, of $J$-holomorphic disks with boundary on the respective Lagrangians and intersecting the complement of the ball. These area-estimates are sharp as can be seen by taking a volume-filling ball embedding into $\CP^2$. Similar statements can be made in higher dimensions.

Recall that we have proved that a split torus of the form $L = T(x,0)$ with $-1 \leqslant x \leqslant 1$ is in neither of the sets $\cl_{\rm Cl}$, $\cl_{\rm Ch}, \cl_{\rm nm}$. 

\begin{question}
Let $L = T(x,0) \subset X_{\alpha}$ be a split torus with $-1 \leqslant x \leqslant 1$. Is $L$ ball-embeddable?
\end{question}

We suspect that the answer is no, but can only prove it in the following special cases. 

\begin{proposition}
Let $L = T(x,0) \subset X_{\alpha}$ be a split torus. If $x=0$ or $\vert x \vert < 1 - \alpha$, then $L$ is not ball-embeddable. 
\end{proposition}

Note that for $\alpha \geqslant 1$, the second part of the statement is empty. This proposition follows from computing displacement energies. On one hand, we have $e(T(0,0))= \infty$ and $e(T(x,0))=1 + \alpha - \vert x \vert$ for $0 < \vert x \vert < 1$. The former follows for example from \cite[Corollary 2.2]{EntPol06}, whereas the latter follows from the proof of Theorem 8.1 in \cite{FOOO13}. On the other hand, $e(\varphi(B^4(r))) = r < 2\alpha$ for all $\alpha < 1$ and all symplectic ball embeddings $\varphi \colon B^4(r) \hookrightarrow X_{\alpha}$. For $\alpha \geqslant 1$, the displacement is less obvious, but also finite. 

The question of ball-embeddability can be asked for \emph{any} Lagrangian torus $L \subset (X,\omega)$. Let us mention a few straightforward results in that direction. In monotone $S^2 \times S^2$, every split torus except for the product of equators is ball-embeddable, as can be seen by taking all possible toric ball embeddings. On the other hand, every non-displaceable torus in monotone $S^2 \times S^2$ is not ball-embeddable. Indeed, every symplectic ball in $S^2 \times S^2$ is displaceable and hence can only contain displaceable tori. Examples of this include the infinitely many monotone tori constructed by Vianna \cite{Via17}, as well as the non-monotone family of non-displaceable tori found by Fukaya--Oh--Ohta--Ono in \cite{FOOO10}. In $\CP^2$ on the other hand, every Lagrangian torus is ball-embeddable. Indeed, every Lagrangian torus in $\CP^2$ can be moved off of the line at infinity by a Hamiltonian isotopy, see \cite[Theorem 1.3]{Wel07} or \cite[Theorem C]{GooIvrRiz16}. For the same reason, every Lagrangian torus in monotone $S^2 \times S^2$ is \emph{polydisk-embeddable}, i.e.\ contained in the image of a symplectic embedding of the product of two disks of equal areas.

\newpage

\bibliographystyle{abbrv}
\bibliography{mybibfile}

\end{document}